\documentclass[twocolumn,english]{IEEEtran}
\usepackage[T1]{fontenc}
\usepackage{color}
\usepackage{babel}
\usepackage{float}
\usepackage{amsthm}
\usepackage{amsmath}
\usepackage{graphicx}
\usepackage[unicode=true,
 bookmarks=true,bookmarksnumbered=true,bookmarksopen=true,bookmarksopenlevel=1,
 breaklinks=false,pdfborder={0 0 0},backref=false,colorlinks=false]
 {hyperref}
\hypersetup{pdftitle={Your Title},
 pdfauthor={Your Name},
 pdfpagelayout=OneColumn, pdfnewwindow=true, pdfstartview=XYZ, plainpages=false}
\usepackage{breakurl}

\makeatletter

\floatstyle{ruled}
\newfloat{algorithm}{tbp}{loa}
\providecommand{\algorithmname}{Algorithm}
\floatname{algorithm}{\protect\algorithmname}

\theoremstyle{plain}
\newtheorem{thm}{\protect\theoremname}
\theoremstyle{plain}
\newtheorem{lem}[thm]{\protect\lemmaname}
 \usepackage{algolyx}
\theoremstyle{plain}
\newtheorem{prop}[thm]{\protect\propositionname}

\ifCLASSOPTIONcompsoc
\usepackage[caption=false,font=normalsize,labelfont=sf,textfont=sf]{subfig}
\else
\usepackage[caption=false,font=footnotesize]{subfig}
\fi

\usepackage{cite}

\keycomment{\#<}{>\#}

\makeatletter
\def\thmhead@plain#1#2#3{%
  \thmname{#1}\thmnumber{\@ifnotempty{#1}{ }\@upn{#2}}%
  \thmnote{ {\the\thm@notefont#3}}}
\let\thmhead\thmhead@plain
\makeatother

\makeatother

\providecommand{\lemmaname}{Lemma}
\providecommand{\propositionname}{Proposition}
\providecommand{\theoremname}{Theorem}

\begin{document}

\title{Sequence Design to Minimize the Weighted Integrated and Peak Sidelobe
Levels}

\author{Junxiao Song, Prabhu Babu, and Daniel P. Palomar, \IEEEmembership{Fellow, IEEE}\thanks{Junxiao Song, Prabhu Babu, and Daniel P. Palomar are with the Hong
Kong University of Science and Technology (HKUST), Hong Kong. E-mail:
\{jsong, eeprabhubabu, palomar\}@ust.hk.}}
\maketitle
\begin{abstract}
Sequences with low aperiodic autocorrelation sidelobes are well known
to have extensive applications in active sensing and communication
systems. In this paper, we consider the problem of minimizing the
weighted integrated sidelobe level (WISL), which can be used to design
sequences with impulse-like autocorrelation and zero (or low) correlation
zone. Two algorithms based on the general majorization-minimization
method are developed to tackle the WISL minimization problem and the
convergence to a stationary point is guaranteed. In addition, the
proposed algorithms can be implemented via fast Fourier transform
(FFT) operations and thus are computationally efficient, and an acceleration
scheme has been considered to further accelerate the algorithms. Moreover,
the proposed methods are extended to optimize the $\ell_{p}$-norm
of the autocorrelation sidelobes, which lead to a way to minimize
the peak sidelobe level (PSL) criterion. Numerical experiments show
that the proposed algorithms can efficiently generate sequences with
virtually zero autocorrelation sidelobes in a specified lag interval
and can also produce very long sequences with much smaller PSL compared
with some well known analytical sequences. \end{abstract}


\section{Introduction}

\IEEEPARstart{S}{equences} with good autocorrelation properties
lie at the heart of many active sensing and communication systems.
Important applications include synchronization of digital communication
systems (e.g., GPS receivers or CDMA cellular systems), pilot sequences
for channel estimation, coded sonar and radar systems and even cryptography
for secure systems \cite{turyn1968sequences,channel_estimation,levanon2004radar,golomb2005signal,he2012waveform}.
In practice, due to the limitations of sequence generation hardware
components (such as the maximum signal amplitude clip of analog-to-digital
converters and power amplifiers), unit-modulus sequences (also known
as polyphase sequences) are of special interest because of their maximum
energy efficiency \cite{he2012waveform}.

Let $\{x_{n}\}_{n=1}^{N}$ denote a complex unit-modulus sequence
of length $N$, then the aperiodic autocorrelations of $\{x_{n}\}_{n=1}^{N}$
are defined as
\begin{equation}
\begin{aligned}r_{k} & =\sum_{n=1}^{N-k}x_{n}^{*}x_{n+k}=r_{-k}^{*},\,k=0,\ldots,N-1.\end{aligned}
\label{eq:aper_r_k}
\end{equation}
The problem of sequence design for good autocorrelation properties
usually arises when small autocorrelation sidelobes (i.e., $k\neq0$)
are required. To measure the goodness of the autocorrelation property
of a sequence, two commonly used metrics are the integrated sidelobe
level (ISL) 
\begin{equation}
{\rm ISL}=\sum_{k=1}^{N-1}|r_{k}|^{2},\label{eq:ISL}
\end{equation}
and the peak sidelobe level (PSL)
\begin{equation}
{\rm PSL}=\max\{|r_{k}|\}_{k=1}^{N-1}.\label{eq:PSL}
\end{equation}

Owing to the practical importance of sequences with low autocorrelation
sidelobes, a lot of effort has been devoted to identifying such sequences.
Binary Barker sequences, with their peak sidelobe level (PSL) no greater
than 1, are perhaps the most well-known such sequences \cite{barker1953group}.
However, it is generally accepted that they do not exist for lengths
greater than 13. In 1965, Golomb and Scholtz \cite{golomb1965generalized}
started to investigate more general sequences called generalized Barker
sequences, which obey the same PSL maximum, but may have complex (polyphase)
elements. Since then, a lot of work has been done to extend the list
of polyphase Barker sequences \cite{zhang1989polyBarker,friese1994polyBarker,brenner1998polyphase,borwein2005polyphase,Polycode2009},
and the longest one ever found is of length 77. It is still unknown
whether there exist longer polyphase Barker sequences. Apart from
searching for longer polyphase Barker sequences, some families of
polyphase sequences with good autocorrelation properties that can
be constructed in closed-form have also been proposed in the literature,
such as the Frank sequences \cite{FrankSequence1963}, the Chu sequences
\cite{chu1972polyphase}, and the Golomb sequences \cite{Golomb1993polyphase}.
It has been shown that the PSL's of these sequences grow almost linearly
with the square root of the length $N$ of the sequences \cite{turyn1967correlation,mow1997aperiodic}. 

In recent years, several optimization based approaches have been proposed
to tackle sequence design problems, see \cite{stoica2009new,ITROX2012,naghsh2013unified,UQP,MISL}.
Among them, \cite{MISL} and \cite{stoica2009new} proposed to design
unit-modulus sequences with low autocorrelation by directly minimizing
the true ISL metric or a simpler criterion that is ``almost equivalent''
to the ISL metric. Efficient algorithms based on fast Fourier transform
(FFT) operations were developed and shown to be capable of producing
very long sequences (of length $10^{4}$ or even larger) with much
lower autocorrelation compared with the Frank sequences and Golomb
sequences. Why the ISL metric was chosen as the objective in the optimization
approaches? It is probably because the ISL metric is more tractable
compared with the PSL metric from an optimization point of view. But
as in the definition of Barker sequences, PSL seems to be the preferred
metric in many cases. So it is of particular interest to also develop
efficient optimization algorithms that can minimize the PSL metric.
Additionally, in some applications, instead of sequences with impulse-like
autocorrelation, zero correlation zone (ZCZ) sequences (with zero
correlations over a smaller range) would suffice \cite{torii2004new,tang2008new}.
In \cite{stoica2009new}, an algorithm named WeCAN (weighted cyclic
algorithm new) was proposed to design sequences with zero or low correlation
zone.

In this paper, we consider the problem of minimizing the weighted
ISL metric, which includes the ISL minimization problem as a special
case and can be used to design zero or low correlation zone sequences
by properly choosing the weights. Two efficient algorithms are developed
based on the general majorization-minimization (MM) method by constructing
two different majorization functions. The proposed algorithms can
be implemented by means of FFT operations and are thus very efficient
in practice. The convergence of the algorithms to a stationary point
is proved. An acceleration scheme is also introduced to further accelerate
the proposed algorithms. We also extend the proposed algorithms to
minimize the $\ell_{p}$-norm of the autocorrelation sidelobes. The
resulting algorithm can be adopted to minimize the PSL metric of unit-modulus
sequences.

The remaining sections of the paper are organized as follows. In Section
\ref{sec:Problem-Formulation}, the problem formulation is presented.
In Sections \ref{sec:Seq-Design-viaMM} and \ref{sec:Diag_majorize},
we first give a brief review of the MM method and then two MM algorithms
are derived, followed by the convergence analysis and an acceleration
scheme in Section \ref{sec:Convergence-Acc}. In Section \ref{sec:Minimizing-Lp-norm},
the algorithms are extended to minimize the $\ell_{p}$-norm of the
autocorrelation sidelobes. Finally, Section \ref{sec:Numerical-Experiments}
presents some numerical results and the conclusions are given in Section
\ref{sec:Conclusion}.

\emph{Notation}: Boldface upper case letters denote matrices, boldface
lower case letters denote column vectors, and italics denote scalars.
$\mathbf{R}$ and $\mathbf{C}$ denote the real field and the complex
field, respectively. $\mathrm{Re}(\cdot)$ and $\mathrm{Im}(\cdot)$
denote the real and imaginary part respectively. ${\rm arg}(\cdot)$
denotes the phase of a complex number. The superscripts $(\cdot)^{T}$,
$(\cdot)^{*}$ and $(\cdot)^{H}$ denote transpose, complex conjugate,
and conjugate transpose, respectively. $\circ$ denotes the Hadamard
product. $X_{i,j}$ denotes the (\emph{i}-th, \emph{j}-th) element
of matrix $\mathbf{X}$ and $x_{i}$ denotes the \emph{i}-th element
of vector $\mathbf{x}$. $\mathbf{X}_{i,:}$ denotes the \emph{i}-th
row of matrix $\mathbf{X}$, $\mathbf{X}_{:,j}$ denotes the \emph{j}-th
column of matrix $\mathbf{X}$, and $\mathbf{X}_{i:j,k:l}$ denotes
the submatrix of $\mathbf{X}$ from $X_{i,k}$ to $X_{j,l}$. $\mathrm{Tr}(\cdot)$
denotes the trace of a matrix. ${\rm diag}(\mathbf{X})$ is a column
vector consisting of all the diagonal elements of $\mathbf{X}$. ${\rm Diag}(\mathbf{x})$
is a diagonal matrix formed with $\mathbf{x}$ as its principal diagonal.
${\rm vec}(\mathbf{X})$ is a column vector consisting of all the
columns of $\mathbf{X}$ stacked. $\mathbf{I}_{n}$ denotes an $n\times n$
identity matrix.

\section{Problem Formulation \label{sec:Problem-Formulation}}

Let $\{x_{n}\}_{n=1}^{N}$ denote the complex unit-modulus sequence
to be designed and $\{r_{k}\}_{k=1}^{N-1}$ be the aperiodic autocorrelations
of $\{x_{n}\}_{n=1}^{N}$ as defined in \eqref{eq:aper_r_k}, then
we define the weighted integrated sidelobe level (WISL) as 
\begin{equation}
{\rm WISL}=\sum_{k=1}^{N-1}w_{k}|r_{k}|^{2},\label{eq:WISL}
\end{equation}
where $w_{k}\geq0,\,k=1,\ldots,N-1.$ It is easy to see that the WISL
metric includes the ISL metric as a special case by simply taking
$w_{k}=1,\,k=1,\ldots,N-1.$ 

The problem of interest in this paper is the following WISL minimization
problem: 
\begin{equation}
\begin{array}{ll}
\underset{x_{n}}{\mathsf{minimize}} & {\rm WISL}\\
\mathsf{subject\;to} & \left|x_{n}\right|=1,\,n=1,\ldots,N,
\end{array}\label{eq:WISL-minimize}
\end{equation}
which includes the ISL minimization problem considered in \cite{MISL}
as a special case and can be used to design zero (or low) correlation
zone sequences by assigning larger weights for the sidelobes that
we want to minimize.

\textcolor{black}{An algorithm named WeCAN was proposed in \cite{stoica2009new}
to tackle this problem. However, instead of directly minimizing the
WISL metric as in \eqref{eq:WISL-minimize}, WeCAN tries to minimize
an ``almost equivalent'' criterion. Moreover, the WeCAN algorithm
requires computing the square-root of an $N\times N$ matrix at the
beginning and $N$ FFT's at each iteration, which could be costly
for large $N$. }In the next section, we will develop algorithms that
directly minimize the original WISL metric, and at the same time are
computationally much more efficient than the WeCAN algorithm.

\section{Sequence Design via Majorization-Minimization\label{sec:Seq-Design-viaMM}}

In this section, we first introduce the general majorization-minimization
(MM) method briefly and then apply it to derive simple algorithms
to solve the problem \eqref{eq:WISL-minimize}.

\subsection{The MM Method\label{sub:MM-Method}}

The MM method refers to the majorization-minimization method, which
is an approach to solve optimization problems that are too difficult
to solve directly. The principle behind the MM method is to transform
a difficult problem into a series of simple problems. Interested readers
may refer to \cite{hunter2004MMtutorial,MM_Stoica} and references
therein for more details (recent generalizations include \cite{razaviyayn2013unified,Aldo2013}).

Suppose we want to minimize $f(\mathbf{x})$ over $\mathcal{X}\subseteq\mathbf{C}^{n}$.
Instead of minimizing the cost function $f(\mathbf{x})$ directly,
the MM approach optimizes a sequence of approximate objective functions
that majorize $f(\mathbf{x})$. More specifically, starting from a
feasible point $\mathbf{x}^{(0)},$ the algorithm produces a sequence
$\{\mathbf{x}^{(k)}\}$ according to the following update rule: 
\begin{equation}
\mathbf{x}^{(k+1)}\in\underset{\mathbf{x}\in\mathcal{X}}{\arg\min}\,\,u(\mathbf{x},\mathbf{x}^{(k)}),\label{eq:major_update}
\end{equation}
where $\mathbf{x}^{(k)}$ is the point generated by the algorithm
at iteration $k,$ and $u(\mathbf{x},\mathbf{x}^{(k)})$ is the majorization
function of $f(\mathbf{x})$ at $\mathbf{x}^{(k)}$. Formally, the
function $u(\mathbf{x},\mathbf{x}^{(k)})$ is said to majorize the
function $f(\mathbf{x})$ at the point $\mathbf{x}^{(k)}$ if 
\begin{eqnarray}
u(\mathbf{x},\mathbf{x}^{(k)}) & \geq & f(\mathbf{x}),\quad\forall\mathbf{x}\in\mathcal{X},\label{eq:major1}\\
u(\mathbf{x}^{(k)},\mathbf{x}^{(k)}) & = & f(\mathbf{x}^{(k)}).\label{eq:major2}
\end{eqnarray}
In other words, function $u(\mathbf{x},\mathbf{x}^{(k)})$ is an upper
bound of $f(\mathbf{x})$ over $\mathcal{X}$ and coincides with $f(\mathbf{x})$
at $\mathbf{x}^{(k)}$.

To summarize, to minimize $f(\mathbf{x})$ over $\mathcal{X}\subseteq\mathbf{C}^{n}$,
the main steps of the majorization-minimization scheme are 
\begin{enumerate}
\item Find a feasible point $\mathbf{x}^{(0)}$ and set $k=0.$ 
\item Construct a function $u(\mathbf{x},\mathbf{x}^{(k)})$ that majorizes
$f(\mathbf{x})$ at $\mathbf{x}^{(k)}$.
\item Let $\mathbf{x}^{(k+1)}\in\underset{\mathbf{x}\in\mathcal{X}}{\arg\min}\,\,u(\mathbf{x},\mathbf{x}^{(k)}).$ 
\item If some convergence criterion is met, exit; otherwise, set $k=k+1$
and go to step (2). 
\end{enumerate}
It is easy to show that with this scheme, the objective value is monotonically
decreasing (nonincreasing) at every iteration, i.e., 
\begin{equation}
f(\mathbf{x}^{(k+1)})\leq u(\mathbf{x}^{(k+1)},\mathbf{x}^{(k)})\leq u(\mathbf{x}^{(k)},\mathbf{x}^{(k)})=f(\mathbf{x}^{(k)}).\label{eq:descent-property}
\end{equation}
The first inequality and the third equality follow from the the properties
of the majorization function, namely \eqref{eq:major1} and \eqref{eq:major2}
respectively and the second inequality follows from \eqref{eq:major_update}.
The monotonicity makes MM algorithms very stable in practice.

\subsection{WISL Minimization via MM\label{sub:WISL-Minimize_sub1}}

To solve the problem \eqref{eq:WISL-minimize} via majorization-minimization,
the key step is to find a majorization function of the objective such
that the majorized problem is easy to solve. For that purpose we first
present a simple result that will be useful later.
\begin{lem}[\cite{MISL}]
\label{lem:majorizer}Let $\mathbf{L}$ be an $n\times n$ Hermitian
matrix and $\mathbf{M}$ be another $n\times n$ Hermitian matrix
such that $\mathbf{M}\succeq\mathbf{L}.$\textup{ }Then for any point
$\mathbf{x}_{0}\in\mathbf{C}^{n}$, the quadratic function $\mathbf{x}^{H}\mathbf{L}\mathbf{x}$
is majorized by $\mathbf{x}^{H}\mathbf{M}\mathbf{x}+2{\rm Re}\left(\mathbf{x}^{H}(\mathbf{L}-\mathbf{M})\mathbf{x}_{0}\right)+\mathbf{x}_{0}^{H}(\mathbf{M}-\mathbf{L})\mathbf{x}_{0}$
at\textup{ $\mathbf{x}_{0}$.}
\end{lem}
Let us define $\mathbf{U}_{k},\,k=0,\ldots,,N-1$ to be $N\times N$
Toeplitz matrices with the $k$th diagonal elements being $1$ and
$0$ elsewhere.Noticing that
\begin{equation}
r_{k}={\rm Tr}(\mathbf{U}_{k}\mathbf{x}\mathbf{x}^{H}),\,k=0,\ldots,N-1,\label{eq:r_k}
\end{equation}
we can rewrite the problem \eqref{eq:WISL-minimize} as 
\begin{equation}
\begin{array}{ll}
\underset{\mathbf{X},\mathbf{x}}{\mathsf{minimize}} & {\displaystyle \sum_{k=1}^{N-1}}w_{k}|{\rm Tr}(\mathbf{U}_{k}\mathbf{X})|^{2}\\
\mathsf{subject\;to} & \mathbf{X}=\mathbf{x}\mathbf{x}^{H}\\
 & \left|x_{n}\right|=1,\,n=1,\ldots,N.
\end{array}\label{eq:WISL-minimize1}
\end{equation}
Let $\mathbf{U}_{-k}=\mathbf{U}_{k}^{T},\,k=1,\ldots,,N-1$, it is
easy to see that ${\rm Tr}(\mathbf{U}_{-k}\mathbf{X})=r_{-k}=r_{k}^{*}$,
so we can rewrite the problem \eqref{eq:WISL-minimize1} in a more
symmetric way:
\begin{equation}
\begin{array}{ll}
\underset{\mathbf{X},\mathbf{x}}{\mathsf{minimize}} & \frac{1}{2}{\displaystyle \sum_{k=1-N}^{N-1}}w_{k}|{\rm Tr}(\mathbf{U}_{k}\mathbf{X})|^{2}\\
\mathsf{subject\;to} & \mathbf{X}=\mathbf{x}\mathbf{x}^{H}\\
 & \left|x_{n}\right|=1,\,n=1,\ldots,N,
\end{array}\label{eq:WISL-min-sym}
\end{equation}
where $w_{-k}=w_{k}$ and $w_{0}=0$. Since ${\rm Tr}(\mathbf{U}_{k}\mathbf{X})={\rm vec}(\mathbf{X})^{H}{\rm vec}(\mathbf{U}_{k})$,
the problem \eqref{eq:WISL-min-sym} can be further rewritten using
${\rm vec}$ notation as (the constant factor $\frac{1}{2}$ is ignored):
\begin{equation}
\begin{array}{ll}
\underset{\mathbf{x},\mathbf{X}}{\mathsf{minimize}} & {\displaystyle \sum_{k=1-N}^{N-1}}w_{k}{\rm vec}(\mathbf{X})^{H}{\rm vec}(\mathbf{U}_{k}){\rm vec}(\mathbf{U}_{k})^{H}{\rm vec}(\mathbf{X})\\
\mathsf{subject\;to} & \mathbf{X}=\mathbf{x}\mathbf{x}^{H}\\
 & \left|x_{n}\right|=1,\,n=1,\ldots,N.
\end{array}\label{eq:WMISL-vecX}
\end{equation}
Let us define 
\begin{equation}
\mathbf{L}={\displaystyle \sum_{k=1-N}^{N-1}}w_{k}{\rm vec}(\mathbf{U}_{k}){\rm vec}(\mathbf{U}_{k})^{H}\label{eq:matrix_L}
\end{equation}
and denote the objective function of \eqref{eq:WMISL-vecX} by $f(\mathbf{X})$,
then $f(\mathbf{X})={\rm vec}(\mathbf{X})^{H}\mathbf{L}{\rm vec}(\mathbf{X})$,
which is clearly a quadratic function in $\mathbf{X}$. According
to Lemma \ref{lem:majorizer}, we can construct a majorization function
of $f(\mathbf{X})$ by simply choosing a matrix $\mathbf{M}$ such
that $\mathbf{M}\succeq\mathbf{L}.$ A simple choice of $\mathbf{M}$
can be $\mathbf{M}=\lambda_{{\rm max}}(\mathbf{L})\mathbf{I}$, where
$\lambda_{{\rm max}}(\mathbf{L})$ is the maximum eigenvalue of $\mathbf{L}$.
Owing to the special structure of $\mathbf{L}$, it can be shown that
$\lambda_{{\rm max}}(\mathbf{L})$ can be computed efficiently in
closed form.
\begin{lem}
\label{lem:lambda_max}Let $\mathbf{L}$ be the matrix defined in
\eqref{eq:matrix_L}. Then the maximum eigenvalue of $\mathbf{L}$
is given by $\lambda_{{\rm max}}(\mathbf{L})=\max_{k}\{w_{k}(N-k)|k=1,\ldots,N-1\}.$\end{lem}
\begin{IEEEproof}
It is easy to see that the set of vectors $\{{\rm vec}(\mathbf{U}_{k})\}_{k=1-N}^{N-1}$
are mutually orthogonal. For $k=1-N,\ldots,N-1,$ we have
\begin{equation}
\begin{aligned}\mathbf{L}{\rm vec}(\mathbf{U}_{k}) & ={\displaystyle \sum_{j=1-N}^{N-1}}w_{j}{\rm vec}(\mathbf{U}_{j}){\rm vec}(\mathbf{U}_{j})^{H}{\rm vec}(\mathbf{U}_{k})\\
 & =w_{k}{\rm vec}(\mathbf{U}_{k}){\rm vec}(\mathbf{U}_{k})^{H}{\rm vec}(\mathbf{U}_{k})\\
 & =w_{k}(N-\left|k\right|){\rm vec}(\mathbf{U}_{k}).
\end{aligned}
\end{equation}
Thus $w_{k}(N-\left|k\right|)$, $k=1-N,\ldots,N-1$ are the nonzero
eigenvalues of $\mathbf{L}$ with corresponding eigenvectors ${\rm vec}(\mathbf{U}_{k}),$
$k=1-N,\ldots,N-1.$ Since $w_{-k}=w_{k}$ and $w_{0}=0$, the maximum
eigenvalue of $\mathbf{L}$ is given by $\max_{k}\{w_{k}(N-k)|k=1,\ldots,N-1\}.$
\end{IEEEproof}
Then given $\mathbf{X}^{(l)}=\mathbf{x}^{(l)}(\mathbf{x}^{(l)})^{H}$
at iteration $l$, by choosing $\mathbf{M}=\lambda_{{\rm max}}(\mathbf{L})\mathbf{I}$
in Lemma \ref{lem:majorizer} we know that the objective of \eqref{eq:WMISL-vecX}
is majorized by the following function at $\mathbf{X}^{(l)}$: 
\begin{equation}
\begin{aligned} & \,\,u_{1}(\mathbf{X},\mathbf{X}^{(l)})\\
= & \,\,\lambda_{{\rm max}}(\mathbf{L}){\rm vec}(\mathbf{X})^{H}{\rm vec}(\mathbf{X})\\
 & +2{\rm Re}\big({\rm vec}(\mathbf{X})^{H}(\mathbf{L}-\lambda_{{\rm max}}(\mathbf{L})\mathbf{I}){\rm vec}(\mathbf{X}^{(l)})\big)\\
 & +{\rm vec}(\mathbf{X}^{(l)})^{H}(\lambda_{{\rm max}}(\mathbf{L})\mathbf{I}-\mathbf{L}){\rm vec}(\mathbf{X}^{(l)}).
\end{aligned}
\label{eq:major_u1}
\end{equation}
Since ${\rm vec}(\mathbf{X})^{H}{\rm vec}(\mathbf{X})=(\mathbf{x}^{H}\mathbf{x})^{2}=N^{2}$,
the first term of \eqref{eq:major_u1} is just a constant. After ignoring
the constant terms, the majorized problem of \eqref{eq:WMISL-vecX}
is given by
\begin{equation}
\begin{array}{ll}
\underset{\mathbf{x},\mathbf{X}}{\mathsf{minimize}} & {\rm Re}\big({\rm vec}(\mathbf{X})^{H}(\mathbf{L}-\lambda_{{\rm max}}(\mathbf{L})\mathbf{I}){\rm vec}(\mathbf{X}^{(l)})\big)\\
\mathsf{subject\;to} & \mathbf{X}=\mathbf{x}\mathbf{x}^{H}\\
 & \left|x_{n}\right|=1,\,n=1,\ldots,N.
\end{array}\label{eq:WMISL-major-prob1}
\end{equation}
Substituting $\mathbf{L}$ in \eqref{eq:matrix_L} back into \eqref{eq:WMISL-major-prob1},
the problems becomes 
\begin{equation}
\begin{array}{ll}
\underset{\mathbf{x},\mathbf{X}}{\mathsf{minimize}} & {\displaystyle \sum_{k=1-N}^{N-1}}w_{k}{\rm Re}({\rm Tr}(\mathbf{U}_{-k}\mathbf{X}^{(l)}){\rm Tr}(\mathbf{U}_{k}\mathbf{X}))\\
 & -\lambda_{{\rm max}}(\mathbf{L}){\rm Tr}(\mathbf{X}^{(l)}\mathbf{X})\\
\mathsf{subject\;to} & \mathbf{X}=\mathbf{x}\mathbf{x}^{H}\\
 & \left|x_{n}\right|=1,\,n=1,\ldots,N,
\end{array}\label{eq:WMISL-major-Tr}
\end{equation}
Since ${\rm Tr}(\mathbf{U}_{-k}\mathbf{X}^{(l)})=r_{-k}^{(l)}$, the
problem \eqref{eq:WMISL-major-Tr} can be rewritten as
\begin{equation}
\begin{array}{ll}
\underset{\mathbf{x},\mathbf{X}}{\mathsf{minimize}} & {\rm Re}\left({\rm Tr}\left({\displaystyle \sum_{k=1-N}^{N-1}}w_{k}r_{-k}^{(l)}\mathbf{U}_{k}\mathbf{X}\right)\right)\\
 & -\lambda_{{\rm max}}(\mathbf{L}){\rm Tr}(\mathbf{X}^{(l)}\mathbf{X})\\
\mathsf{subject\;to} & \mathbf{X}=\mathbf{x}\mathbf{x}^{H}\\
 & \left|x_{n}\right|=1,\,n=1,\ldots,N,
\end{array}\label{eq:WMISL-major-Tr-r}
\end{equation}
which can be further simplified as 
\begin{equation}
\begin{array}{ll}
\underset{\mathbf{x}}{\mathsf{minimize}} & \mathbf{x}^{H}\left(\mathbf{R}-\lambda_{{\rm max}}(\mathbf{L})\mathbf{x}^{(l)}(\mathbf{x}^{(l)})^{H}\right)\mathbf{x}\\
\mathsf{subject\;to} & \left|x_{n}\right|=1,\,n=1,\ldots,N,
\end{array}\label{eq:toep}
\end{equation}
where
\begin{equation}
\begin{aligned}\mathbf{R} & ={\displaystyle \sum_{k=1-N}^{N-1}}w_{k}r_{-k}^{(l)}\mathbf{U}_{k}\\
 & =\begin{bmatrix}0 & w_{1}r_{-1}^{(l)} & \ldots & w_{N-1}r_{1-N}^{(l)}\\
w_{1}r_{1}^{(l)} & 0 & \ddots & \vdots\\
\vdots & \ddots & \ddots & w_{1}r_{-1}^{(l)}\\
w_{N-1}r_{N-1}^{(l)} & \ldots & w_{1}r_{1}^{(l)} & 0
\end{bmatrix}
\end{aligned}
\label{eq:R_mat}
\end{equation}
is a Hermitian Toeplitz matrix, and $\{r_{k}^{(l)}\}_{k=1-N}^{N-1}$
are the autocorrelations of the sequence $\{x_{n}^{(l)}\}_{n=1}^{N}$.

It is clear that the objective function in \eqref{eq:toep} is quadratic
in $\mathbf{x}$, but the problem \eqref{eq:toep} is still hard to
solve directly. So we propose to majorize the objective function of
problem \eqref{eq:toep} at $\mathbf{x}^{(l)}$ again to further simplify
the problem that we need to solve at each iteration. Similarly, to
construct a majorization function of the objective, we need to find
a matrix $\mathbf{M}$ such that $\mathbf{M}\succeq\mathbf{R}-\lambda_{{\rm max}}(\mathbf{L})\mathbf{x}^{(l)}(\mathbf{x}^{(l)})^{H}.$
As one choice, one may choose $\mathbf{M}=\lambda_{{\rm max}}\left(\mathbf{R}-\lambda_{{\rm max}}(\mathbf{L})\mathbf{x}^{(l)}(\mathbf{x}^{(l)})^{H}\right)\mathbf{I}$,
as in the first majorization step. But in this case, to compute the
maximum eigenvalue of the matrix $\mathbf{R}-\lambda_{{\rm max}}(\mathbf{L})\mathbf{x}^{(l)}(\mathbf{x}^{(l)})^{H}$,
some iterative algorithms are needed, in contrast to the simple closed
form expression in the first majorization step. To maintain the simplicity
and the computational efficiency of the algorithm, here we propose
to use some upper bound of $\lambda_{{\rm max}}\left(\mathbf{R}-\lambda_{{\rm max}}(\mathbf{L})\mathbf{x}^{(l)}(\mathbf{x}^{(l)})^{H}\right)$
that can be easily computed instead. To derive the upper bound, we
first introduce a useful result regarding the bounds of the extreme
eigenvalues of Hermitian Toeplitz matrices \cite{eig_localization}.
\begin{lem}
\label{lem:eig_bounds}Let $\mathbf{T}$ be an $N\times N$ Hermitian
Toeplitz matrix defined by $\{t_{k}\}_{k=0}^{N-1}$ as follows
\[
\mathbf{T}=\begin{bmatrix}t_{0} & t_{1}^{*} & \ldots & t_{N-1}^{*}\\
t_{1} & t_{0} & \ddots & \vdots\\
\vdots & \ddots & \ddots & t_{1}^{*}\\
t_{N-1} & \ldots & t_{1} & t_{0}
\end{bmatrix}
\]
and $\mathbf{F}$ be a $2N\times2N$ FFT matrix with $F_{m,n}=e^{-j\frac{2mn\pi}{2N}},0\leq m,n<2N$.
Let $\mathbf{c}=[t_{0},t_{1},\cdots,t_{N-1},0,t_{N-1}^{*},\cdots,t_{1}^{*}]^{T}$
and $\boldsymbol{\mu}=\mathbf{F}\mathbf{c}$ be the discrete Fourier
transform of $\mathbf{c}$. Then 
\begin{eqnarray}
\lambda_{{\rm max}}(\mathbf{T}) & \leq & \frac{1}{2}\left(\max_{1\leq i\leq N}\mu_{2i}+\max_{1\leq i\leq N}\mu_{2i-1}\right),\label{eq:eig_up}\\
\lambda_{{\rm min}}(\mathbf{T}) & \geq & \frac{1}{2}\left(\min_{1\leq i\leq N}\mu_{2i}+\min_{1\leq i\leq N}\mu_{2i-1}\right).\label{eq:eig_low}
\end{eqnarray}
\end{lem}
\begin{IEEEproof}
See \cite{eig_localization}.
\end{IEEEproof}
Since the matrix $\mathbf{R}$ is Hermitian Toeplitz, according to
Lemma \ref{lem:eig_bounds}, we know that 
\begin{equation}
\lambda_{{\rm max}}\left(\mathbf{R}\right)\leq\frac{1}{2}\left(\max_{1\leq i\leq N}\mu_{2i}+\max_{1\leq i\leq N}\mu_{2i-1}\right),\label{eq:eig_upper_bound_R}
\end{equation}
where $\boldsymbol{\mu}=\mathbf{F}\mathbf{c}$ and 
\begin{equation}
\mathbf{c}=[0,w_{1}r_{1}^{(l)},\ldots,w_{N-1}r_{N-1}^{(l)},0,w_{N-1}r_{1-N}^{(l)},\ldots,w_{1}r_{-1}^{(l)}]^{T}.\label{eq:c_coef}
\end{equation}
Let us denote the upper bound of $\lambda_{{\rm max}}\left(\mathbf{R}\right)$
at the right hand side of \eqref{eq:eig_upper_bound_R} by $\lambda_{u}$,
i.e., 
\begin{equation}
\lambda_{u}=\frac{1}{2}\left(\max_{1\leq i\leq N}\mu_{2i}+\max_{1\leq i\leq N}\mu_{2i-1}\right).\label{eq:lambda_u}
\end{equation}
 Since $\lambda_{{\rm max}}(\mathbf{L})\geq0$, it is easy to see
that 
\begin{equation}
\lambda_{u}\geq\lambda_{{\rm max}}\left(\mathbf{R}\right)\geq\lambda_{{\rm max}}\left(\mathbf{R}-\lambda_{{\rm max}}(\mathbf{L})\mathbf{x}^{(l)}(\mathbf{x}^{(l)})^{H}\right).
\end{equation}
Thus, we may choose $\mathbf{M}=\lambda_{u}\mathbf{I}$ in Lemma \ref{lem:majorizer}
and the objective of \eqref{eq:toep} is majorized by 
\begin{equation}
\begin{aligned} & \,\,u_{2}(\mathbf{x},\mathbf{x}^{(l)})\\
= & \,\,\lambda_{u}\mathbf{x}^{H}\mathbf{x}\!+\!2{\rm Re}\left(\mathbf{x}^{H}(\mathbf{R}\!-\!\lambda_{{\rm max}}(\mathbf{L})\mathbf{x}^{(l)}(\mathbf{x}^{(l)})^{H}\!-\!\lambda_{u}\mathbf{I})\mathbf{x}^{(l)}\right)\\
 & +(\mathbf{x}^{(l)})^{H}(\lambda_{u}\mathbf{I}-\mathbf{R}\!+\!\lambda_{{\rm max}}(\mathbf{L})\mathbf{x}^{(l)}(\mathbf{x}^{(l)})^{H})\mathbf{x}^{(l)}.
\end{aligned}
\label{eq:major_u2}
\end{equation}

Since $\mathbf{x}^{H}\mathbf{x}=N,$ the first term of \eqref{eq:major_u2}
is a constant. Again by ignoring the constant terms, we have the majorized
problem of \eqref{eq:toep}:
\begin{equation}
\begin{array}{ll}
\underset{\mathbf{x}}{\mathsf{minimize}} & {\rm Re}\left(\mathbf{x}^{H}(\mathbf{R}-\lambda_{{\rm max}}(\mathbf{L})\mathbf{x}^{(l)}(\mathbf{x}^{(l)})^{H}-\lambda_{u}\mathbf{I})\mathbf{x}^{(l)}\right)\\
\mathsf{subject\;to} & \left|x_{n}\right|=1,\,n=1,\ldots,N,
\end{array}\label{eq:linear_prob}
\end{equation}
which can be rewritten as
\begin{equation}
\begin{array}{ll}
\underset{\mathbf{x}}{\mathsf{minimize}} & \left\Vert \mathbf{x}-\mathbf{y}\right\Vert _{2}\\
\mathsf{subject\;to} & \left|x_{n}\right|=1,\,n=1,\ldots,N,
\end{array}\label{eq:toep3}
\end{equation}
where
\begin{eqnarray}
\mathbf{y} & = & -(\mathbf{R}-\lambda_{{\rm max}}(\mathbf{L})\mathbf{x}^{(l)}(\mathbf{x}^{(l)})^{H}-\lambda_{u}\mathbf{I})\mathbf{x}^{(l)}\nonumber \\
 & = & \left(\lambda_{{\rm max}}(\mathbf{L})N+\lambda_{u}\right)\mathbf{x}^{(l)}-\mathbf{R}\mathbf{x}^{(l)}.\label{eq:y_k2}
\end{eqnarray}
It is easy to see that the problem \eqref{eq:toep3} has a closed
form solution, which is given by 
\begin{equation}
x_{n}=e^{j{\rm arg}(y_{n})},\,n=1,\ldots,N.\label{eq:WISL_minimizer}
\end{equation}

Note that although we have applied the majorization-minimization scheme
twice at the point $\mathbf{x}^{(l)}$, it can be viewed as directly
majorizing the objective function of \eqref{eq:WISL-minimize} at
$\mathbf{x}^{(l)}$ by the following function:
\begin{equation}
\begin{aligned} & \,\,u(\mathbf{x},\mathbf{x}^{(l)})\\
= & \,\,u_{2}(\mathbf{x},\mathbf{x}^{(l)})\!+\!\lambda_{{\rm max}}(\mathbf{L})N^{2}\!-\!{\displaystyle \sum_{k=1}^{N-1}}w_{k}|r_{k}^{(l)}|^{2}\\
= & \,\,-2{\rm Re}\left(\mathbf{x}^{H}\mathbf{y}\right)\!+\!2N(\lambda_{{\rm max}}(\mathbf{L})N+\lambda_{u})\!-\!3{\displaystyle \sum_{k=1}^{N-1}}w_{k}|r_{k}^{(l)}|^{2},
\end{aligned}
\label{eq:upper_u}
\end{equation}
and the minimizer of $u(\mathbf{x},\mathbf{x}^{(l)})$ over the constraint
set is given by \eqref{eq:WISL_minimizer}. 

According to the steps of the majorization-minimization scheme described
in section \ref{sub:MM-Method}, we can now readily have a straightforward
implementation of the algorithm, which at each iteration computes
$\mathbf{y}$ according to \eqref{eq:y_k2} and update $\mathbf{x}$
via \eqref{eq:WISL_minimizer}. It is easy to see that the main cost
is the computation of $\mathbf{y}$ in \eqref{eq:y_k2}. To obtain
an efficient implementation of the algorithm, here we further explore
the special structures of the matrices involved in the computation
of $\mathbf{y}$. 

We first notice that to compute $\lambda_{u},$ we need to compute
the FFT of the vector $\mathbf{c}$ in \eqref{eq:c_coef} and the
autocorrelations $\{r_{k}^{(l)}\}_{k=1-N}^{N-1}$ of $\{x_{n}^{(l)}\}_{n=1}^{N}$
are needed to form the vector $\mathbf{c}$. It is well known that
the autocorrelations can be computed efficiently via FFT (IFFT) operations,
i.e.,
\begin{equation}
\begin{aligned} & \,\,[r_{0}^{(l)},r_{1}^{(l)},\ldots,r_{N-1}^{(l)},0,r_{1-N}^{(l)},\ldots,r_{-1}^{(l)}]^{T}\\
= & \,\,\frac{1}{2N}\mathbf{F}^{H}\left|\mathbf{F}[\mathbf{x}^{(l)T},\mathbf{0}_{1\times N}]^{T}\right|^{2},
\end{aligned}
\label{eq:auto_FFT}
\end{equation}
where $\mathbf{F}$ is the $2N\times2N$ FFT matrix and $\left|\cdot\right|^{2}$
denotes the element-wise absolute-squared value. Next we present
another simple result regarding Hermitian Toeplitz matrices that can
be used to compute the matrix vector multiplication $\mathbf{R}\mathbf{x}^{(l)}$
efficiently via FFT (IFFT). 
\begin{lem}
\label{lem:diagonal}Let $\mathbf{T}$ be an $N\times N$ Hermitian
Toeplitz matrix defined as follows
\[
\mathbf{T}=\begin{bmatrix}t_{0} & t_{1}^{*} & \ldots & t_{N-1}^{*}\\
t_{1} & t_{0} & \ddots & \vdots\\
\vdots & \ddots & \ddots & t_{1}^{*}\\
t_{N-1} & \ldots & t_{1} & t_{0}
\end{bmatrix}
\]
and $\mathbf{F}$ be a $2N\times2N$ FFT matrix with $F_{m,n}=e^{-j\frac{2mn\pi}{2N}},0\leq m,n<2N$.
Then $\mathbf{T}$ can be decomposed as $\mathbf{T}=\frac{1}{2N}\mathbf{F}_{:,1:N}^{H}{\rm Diag}(\mathbf{F}\mathbf{c})\mathbf{F}_{:,1:N}$,
where $\mathbf{c}=[t_{0},t_{1},\cdots,t_{N-1},0,t_{N-1}^{*},\cdots,t_{1}^{*}]^{T}$.\end{lem}
\begin{IEEEproof}
The $N\times N$ Hermitian Toeplitz matrix $\mathbf{T}$ can be embedded
in a circulant matrix $\mathbf{C}$ of dimension $2N\times2N$ as
follows:
\begin{equation}
\mathbf{C}=\left[\begin{array}{cc}
\mathbf{T} & \mathbf{W}\\
\mathbf{W} & \mathbf{T}
\end{array}\right],
\end{equation}
where
\begin{equation}
\mathbf{W}=\begin{bmatrix}0 & t_{N-1} & \cdots & t_{1}\\
t_{N-1}^{*} & 0 & \ddots & \vdots\\
\vdots & \ddots & \ddots & t_{N-1}\\
t_{1}^{*} & \cdots & t_{N-1}^{*} & 0
\end{bmatrix}.
\end{equation}
The circulant matrix $\mathbf{C}$ can be diagonalized by the FFT
matrix \cite{Gray2006}, i.e.,
\begin{equation}
\mathbf{C}=\frac{1}{2N}\mathbf{F}^{H}{\rm Diag}(\mathbf{F}\mathbf{c})\mathbf{F},\label{eq:circ_decomp}
\end{equation}
where $\mathbf{c}$ is the first column of $\mathbf{C},$ i.e., $\mathbf{c}=[t_{0},t_{1},\cdots,t_{N-1},0,t_{N-1}^{*},\cdots,t_{1}^{*}]^{T}.$
Since the matrix $\mathbf{T}$ is just the upper left $N\times N$
block of $\mathbf{C}$, we can easily obtain $\mathbf{T}=\frac{1}{2N}\mathbf{F}_{:,1:N}^{H}{\rm Diag}(\mathbf{F}\mathbf{c})\mathbf{F}_{:,1:N}$.
\end{IEEEproof}
Since the matrix $\mathbf{R}$ is Hermitian Toeplitz, from Lemma \ref{lem:diagonal}
we easily have 
\begin{equation}
\mathbf{R}=\frac{1}{2N}\mathbf{F}_{:,1:N}^{H}{\rm Diag}(\mathbf{F}\mathbf{c})\mathbf{F}_{:,1:N},\label{eq:R_decomp}
\end{equation}
where $\mathbf{c}$ is the same as the one defined in \eqref{eq:c_coef},
which can be reused here. With the decomposition of $\mathbf{R}$
given in \eqref{eq:R_decomp}, it is easy to see that the matrix vector
multiplication $\mathbf{R}\mathbf{x}^{(l)}$ can be performed by means
of FFT (IFFT) operations.

Now we are ready to summarize the overall algorithm and it is given
in Algorithm \ref{alg:MWISL}. The algorithm falls into the general
framework of MM algorithms and thus preserves the monotonicity of
such algorithms. It is also worth noting that the per iteration computation
of the algorithm is dominated by four FFT (IFFT) operations and thus
computationally very efficient.

\begin{algorithm}[tbh]
\begin{algor}[1]
\item [{Require:}] \begin{raggedright}
sequence length $N$, weights $\{w_{k}\geq0\}_{k=1}^{N-1}$
\par\end{raggedright}
\item [{{*}}] Set $l=0$, initialize $\mathbf{x}^{(0)}$. 
\item [{{*}}] $\lambda_{L}=\max_{k}\{w_{k}(N-k)|k=1,\ldots,N-1\}$
\item [{repeat}]~

\begin{algor}[1]
\item [{{*}}] $\mathbf{f}=\mathbf{F}[\mathbf{x}^{(l)T},\mathbf{0}_{1\times N}]^{T}$
\item [{{*}}] $\mathbf{r}=\frac{1}{2N}\mathbf{F}^{H}\left|\mathbf{f}\right|^{2}$
\item [{{*}}] $\mathbf{c}=\mathbf{r}\circ[0,w_{1},\ldots,w_{N-1},0,w_{N-1},\ldots,w_{1}]^{T}$
\item [{{*}}] $\boldsymbol{\mu}=\mathbf{F}\mathbf{c}$
\item [{{*}}] $\lambda_{u}=\frac{1}{2}\big(\underset{1\leq i\leq N}{\max}\mu_{2i}+\underset{1\leq i\leq N}{\max}\mu_{2i-1}\big)$
\item [{{*}}] $\mathbf{y}=\mathbf{x}^{(l)}-\frac{\mathbf{F}_{:,1:N}^{H}(\boldsymbol{\mu}\circ\mathbf{f})}{2N(\lambda_{L}N+\lambda_{u})}$
\item [{{*}}] $x_{n}^{(l+1)}=e^{j{\rm arg}(y_{n})},\,n=1,\ldots,N$ 
\item [{{*}}] $l\leftarrow l+1$ 
\end{algor}
\item [{until}] convergence
\end{algor}
\protect\caption{\label{alg:MWISL}MWISL - Monotonic minimizer for Weighted ISL.}
\end{algorithm}

\section{WISL Minimization with an Improved Majorization Function\label{sec:Diag_majorize}}

As described in the previous section, the proposed algorithm is based
on the majorization-minimization principle, and the nature of the
majorization functions usually dictate the performance of the algorithm.
In this section, we further explore the special structure of the problem
and construct a different majorization function.

Notice that to obtain the simple algorithm in the previous section,
a key point is that the first term of the majorization function $u_{1}(\mathbf{X},\mathbf{X}^{(l)})$
in \eqref{eq:major_u1} is just a constant and can be ignored, which
removes the higher order term in the objective of the majorized problem
\eqref{eq:WMISL-major-prob1}. In the previous section, we have chosen
$\mathbf{M}=\lambda_{{\rm max}}(\mathbf{L})\mathbf{I}\succeq\mathbf{L}$
such that ${\rm vec}(\mathbf{X})^{H}\mathbf{M}{\rm vec}(\mathbf{X})=\lambda_{{\rm max}}(\mathbf{L}){\rm vec}(\mathbf{X})^{H}{\rm vec}(\mathbf{X})$
is a constant over the constraint set. But it is easy to see that,
to ensure the term ${\rm vec}(\mathbf{X})^{H}\mathbf{M}{\rm vec}(\mathbf{X})$
being constant, it is sufficient to require $\mathbf{M}$ being diagonal,
i.e., choosing $\mathbf{M}={\rm Diag}(\mathbf{b})\succeq\mathbf{L}$.
To construct a tight majorization function, here we consider the following
problem
\begin{equation}
\begin{array}{ll}
\underset{\mathbf{b}}{\mathsf{minimize}} & {\rm Tr}({\rm Diag}(\mathbf{b})-\mathbf{L})\\
\mathsf{subject\;to} & {\rm Diag}(\mathbf{b})\succeq\mathbf{L},
\end{array}\label{eq:diag_SDP_L}
\end{equation}
i.e., we choose the diagonal matrix ${\rm Diag}(\mathbf{b})$ that
minimizes the sum of eigenvalues of the difference ${\rm Diag}(\mathbf{b})-\mathbf{L}.$
Since $\mathbf{L}$ is a constant matrix, the problem \eqref{eq:diag_SDP_L}
can be written as 
\begin{equation}
\begin{array}{ll}
\underset{\mathbf{b}}{\mathsf{minimize}} & \mathbf{b}^{T}\mathbf{1}_{n}\\
\mathsf{subject\;to} & {\rm Diag}(\mathbf{b})\succeq\mathbf{L}.
\end{array}\label{eq:diag_SDP_L2}
\end{equation}
The problem \eqref{eq:diag_SDP_L2} is an SDP (semidefinite programming),
and there is no closed form solution in general. However, due to the
special properties (i.e., symmetry and nonnegativity) of the matrix
$\mathbf{L}$ in \eqref{eq:matrix_L}, it can be shown that the problem
\eqref{eq:diag_SDP_L2} can be solved in closed form and it is given
in the following lemma.
\begin{lem}
\label{lem:diag_SDP}Let $\mathbf{L}$ be an $n\times n$ real symmetric
nonnegative matrix. Then the problem
\begin{equation}
\begin{array}{ll}
\underset{\mathbf{b}}{\mathsf{minimize}} & \mathbf{b}^{T}\mathbf{1}_{n}\\
\mathsf{subject\;to} & {\rm Diag}(\mathbf{b})\succeq\mathbf{L}
\end{array}\label{eq:diag_SDP_prob}
\end{equation}
admits the following optimal solution:
\begin{equation}
\mathbf{b}^{\star}=\mathbf{L}\mathbf{1}_{n}.
\end{equation}
\end{lem}
\begin{IEEEproof}
Since ${\rm Diag}(\mathbf{b})-\mathbf{L}\succeq\mathbf{0}$, we have
\[
\mathbf{x}^{T}\left({\rm Diag}(\mathbf{b})-\mathbf{L}\right)\mathbf{x}\geq0,\,\forall\mathbf{x}.
\]
By choosing $\mathbf{x}=\mathbf{1}_{n}$, we get $\mathbf{b}^{T}\mathbf{1}_{n}\geq\mathbf{1}_{n}^{T}\mathbf{L}\mathbf{1}_{n}$,
with the equality achieved by $\mathbf{b}=\mathbf{L}\mathbf{1}_{n}$.
So it remains to show that $\mathbf{b}=\mathbf{L}\mathbf{1}_{n}$
is feasible, i.e., ${\rm Diag}(\mathbf{L}\mathbf{1}_{n})\succeq\mathbf{L}$.
Since for any $\mathbf{x}\in\mathbf{R}^{n}$, we have 
\begin{equation}
\begin{aligned} & \,\,\mathbf{x}^{T}\left({\rm Diag}(\mathbf{L}\mathbf{1}_{n})-\mathbf{L}\right)\mathbf{x}\\
= & \,\,\sum_{i}\big(x_{i}^{2}\sum_{j}L_{i,j}\big)-\sum_{i,j}x_{i}L_{i,j}x_{j}\\
= & \,\,\sum_{i,j}\big(L_{i,j}x_{i}^{2}-x_{i}L_{i,j}x_{j}\big)\\
= & \,\,\sum_{i,j}\frac{1}{2}L_{i,j}\big(x_{i}^{2}+x_{j}^{2}-2x_{i}x_{j}\big)\\
= & \,\,\sum_{i,j}\frac{1}{2}L_{i,j}\big(x_{i}-x_{j}\big)^{2}\\
\geq & 0,
\end{aligned}
\end{equation}
where the third equality follows from the symmetry of $\mathbf{L}$
and the last inequality follows from the fact that $\mathbf{L}$ is
nonnegative, the proof is complete. 
\end{IEEEproof}
Then given $\mathbf{X}^{(l)}=\mathbf{x}^{(l)}(\mathbf{x}^{(l)})^{H}$
at iteration $l$, by choosing $\mathbf{M}={\rm Diag}(\mathbf{b})={\rm Diag}(\mathbf{L}\mathbf{1})$
in Lemma \ref{lem:majorizer}, the objective of \eqref{eq:WMISL-vecX}
is majorized by the following function at $\mathbf{X}^{(l)}$:
\begin{equation}
\begin{aligned} & \,\,\tilde{u}_{1}(\mathbf{X},\mathbf{X}^{(l)})\\
= & \,\,{\rm vec}(\mathbf{X})^{H}{\rm Diag}(\mathbf{b}){\rm vec}(\mathbf{X})\\
 & +2{\rm Re}\big({\rm vec}(\mathbf{X})^{H}(\mathbf{L}-{\rm Diag}(\mathbf{b})){\rm vec}(\mathbf{X}^{(l)})\big)\\
 & +{\rm vec}(\mathbf{X}^{(l)})^{H}({\rm Diag}(\mathbf{b})-\mathbf{L}){\rm vec}(\mathbf{X}^{(l)}).
\end{aligned}
\label{eq:major_u1-diag}
\end{equation}
As discussed above, the first term of \eqref{eq:major_u1-diag} is
a constant and by ignoring the constant terms we have the majorized
problem given as follows:
\begin{equation}
\begin{array}{ll}
\underset{\mathbf{x},\mathbf{X}}{\mathsf{minimize}} & {\rm Re}\big({\rm vec}(\mathbf{X})^{H}(\mathbf{L}-{\rm Diag}(\mathbf{b})){\rm vec}(\mathbf{X}^{(l)})\big)\\
\mathsf{subject\;to} & \mathbf{X}=\mathbf{x}\mathbf{x}^{H}\\
 & \left|x_{n}\right|=1,\,n=1,\ldots,N.
\end{array}\label{eq:major_first_diag}
\end{equation}
Since
\begin{equation}
\begin{aligned} & \,\,{\rm Re}\big({\rm vec}(\mathbf{X})^{H}{\rm Diag}(\mathbf{b}){\rm vec}(\mathbf{X}^{(l)})\big)\\
= & \,\,{\rm Re}\big({\rm vec}(\mathbf{X})^{H}\big(\mathbf{b}\circ{\rm vec}(\mathbf{X}^{(l)})\big)\big)\\
= & \,\,{\rm Tr}\big(\mathbf{x}\mathbf{x}^{H}{\rm mat}\big(\mathbf{b}\circ{\rm vec}(\mathbf{x}^{(l)}(\mathbf{x}^{(l)})^{H})\big)\big)\\
= & \,\,\mathbf{x}^{H}\left({\rm mat}(\mathbf{b})\circ\big(\mathbf{x}^{(l)}(\mathbf{x}^{(l)})^{H}\big)\right)\mathbf{x},
\end{aligned}
\label{eq:u1-diag-sec}
\end{equation}
where ${\rm mat}(\cdot)$ is the inverse operation of ${\rm vec}(\cdot),$
following the derivation in the previous section, we can simplify
\eqref{eq:major_first_diag} to
\begin{equation}
\begin{array}{ll}
\underset{\mathbf{x}}{\mathsf{minimize}} & \mathbf{x}^{H}\left(\mathbf{R}-\mathbf{B}\circ\big(\mathbf{x}^{(l)}(\mathbf{x}^{(l)})^{H}\big)\right)\mathbf{x}\\
\mathsf{subject\;to} & \left|x_{n}\right|=1,\,n=1,\ldots,N,
\end{array}\label{eq:toep-diag}
\end{equation}
where $\mathbf{R}$ is defined in \eqref{eq:R_mat} and
\[
\begin{aligned}\mathbf{B} & ={\rm mat}(\mathbf{b})\\
 & ={\rm mat}(\mathbf{L}\mathbf{1})\\
 & ={\rm mat}\left({\displaystyle \sum_{k=1-N}^{N-1}}w_{k}{\rm vec}(\mathbf{U}_{k}){\rm vec}(\mathbf{U}_{k})^{H}\mathbf{1}\right)\\
 & ={\displaystyle \sum_{k=1-N}^{N-1}}w_{k}(N-\left|k\right|)\mathbf{U}_{k}.
\end{aligned}
\]
The objective in \eqref{eq:toep-diag} is quadratic in $\mathbf{x}$
and similar as before we would like to find a matrix $\mathbf{M}$
such that $\mathbf{M}\succeq\mathbf{R}-\mathbf{B}\circ\big(\mathbf{x}^{(l)}(\mathbf{x}^{(l)})^{H}\big).$
For the same reason as in the previous section, here we use some upper
bound of $\lambda_{{\rm max}}\left(\mathbf{R}-\mathbf{B}\circ\big(\mathbf{x}^{(l)}(\mathbf{x}^{(l)})^{H}\big)\right)$
to construct the matrix $\mathbf{M}$. We first present two results
that will be used to derive the upper bound. The first result indicates
some relations between Hadamard and conventional matrix multiplication
\cite{roger1994topics}.
\begin{lem}
\label{lem:Hadamard}Let $\mathbf{A}$, $\mathbf{B}\in\mathbf{C}^{m\times n}$
and let $\mathbf{x}\in\mathbf{C}^{n}$. Then the $i$th diagonal entry
of the matrix $\mathbf{A}{\rm Diag}(\mathbf{x})\mathbf{B}^{T}$ coincides
with the $i$th entry of the vector $\big(\mathbf{A}\circ\mathbf{B}\big)\mathbf{x}$,
$i=1,\ldots,m$, i.e., 
\begin{equation}
\big(\mathbf{A}\circ\mathbf{B}\big)\mathbf{x}={\rm diag}\left(\mathbf{A}{\rm Diag}(\mathbf{x})\mathbf{B}^{T}\right).
\end{equation}

\end{lem}
Then the second result follows, which reveals a fact regarding the
eigenvalues of the matrix $\mathbf{B}\circ\big(\mathbf{x}^{(l)}(\mathbf{x}^{(l)})^{H}\big)$.
\begin{lem}
\label{lem:eig_set}Let $\mathbf{B}$ be an $N\times N$ matrix and
$\mathbf{x}\in\mathbf{C}^{N}$ with $\left|x_{n}\right|=1,\,n=1,\ldots,N$.
Then $\mathbf{B}\circ(\mathbf{x}\mathbf{x}^{H})$ and $\mathbf{B}$
share the same set of eigenvalues. \end{lem}
\begin{IEEEproof}
Suppose $\lambda$ is an eigenvalue of $\mathbf{B}$ and $\mathbf{z}$
is the corresponding eigenvector, i.e., $\mathbf{B}\mathbf{z}=\lambda\mathbf{z}$,
then 
\begin{equation}
\begin{aligned} & \,\,\big(\mathbf{B}\circ(\mathbf{x}\mathbf{x}^{H})\big)(\mathbf{x}\circ\mathbf{z})\\
= & \,\,{\rm diag}\big(\mathbf{B}{\rm Diag}(\mathbf{x}\circ\mathbf{z})(\mathbf{x}\mathbf{x}^{H})^{T}\big)\\
= & \,\,{\rm diag}\big(\mathbf{B}(\mathbf{x}\circ\mathbf{z}\circ\mathbf{x}^{*})\mathbf{x}^{T}\big)\\
= & \,\,{\rm diag}\big(\mathbf{B}\mathbf{z}\mathbf{x}^{T}\big)\\
= & \,\,{\rm diag}\big(\lambda\mathbf{z}\mathbf{x}^{T}\big)\\
= & \,\,\lambda(\mathbf{x}\circ\mathbf{z}),
\end{aligned}
\label{eq:same_eig}
\end{equation}
which means $\lambda$ is also an eigenvalue of the matrix $\mathbf{B}\circ(\mathbf{x}\mathbf{x}^{H})$,
with the corresponding eigenvector given by $(\mathbf{x}\circ\mathbf{z})$. 
\end{IEEEproof}
With Lemma \ref{lem:eig_set}, we have 
\begin{equation}
\begin{aligned} & \,\,\lambda_{{\rm max}}\left(\mathbf{R}-\mathbf{B}\circ\big(\mathbf{x}^{(l)}(\mathbf{x}^{(l)})^{H}\big)\right)\\
\leq & \,\,\lambda_{{\rm max}}(\mathbf{R})-\lambda_{{\rm min}}\left(\mathbf{B}\circ\big(\mathbf{x}^{(l)}(\mathbf{x}^{(l)})^{H}\big)\right)\\
= & \,\,\lambda_{{\rm max}}(\mathbf{R})-\lambda_{{\rm min}}\left(\mathbf{B}\right).
\end{aligned}
\label{eq:eig_u_diag}
\end{equation}
Noticing that the matrix $\mathbf{B}$ is symmetric Toeplitz, according
to Lemma \ref{lem:eig_bounds}, we know that 
\begin{equation}
\lambda_{{\rm min}}(\mathbf{B})\geq\frac{1}{2}\left(\min_{1\leq i\leq N}\nu_{2i}+\min_{1\leq i\leq N}\nu_{2i-1}\right),
\end{equation}
where $\boldsymbol{\nu}=\mathbf{F}\tilde{\mathbf{w}}$ and 
\begin{equation}
\tilde{\mathbf{w}}=[0,w_{1}(N-1),\ldots,w_{N-1},0,w_{N-1},\ldots,w_{1}(N-1)]^{T}.
\end{equation}
By defining 
\begin{equation}
\lambda_{B}=\frac{1}{2}\left(\min_{1\leq i\leq N}\nu_{2i}+\min_{1\leq i\leq N}\nu_{2i-1}\right),
\end{equation}
we now have 
\begin{equation}
\lambda_{{\rm max}}\left(\mathbf{R}-\mathbf{B}\circ\big(\mathbf{x}^{(l)}(\mathbf{x}^{(l)})^{H}\big)\right)\leq\lambda_{u}-\lambda_{B},
\end{equation}
where $\lambda_{u}$ is defined in \eqref{eq:lambda_u} and by choosing
$\mathbf{M}=(\lambda_{u}-\lambda_{B})\mathbf{I}$ in Lemma \ref{lem:majorizer}
we know that the objective of \eqref{eq:toep-diag} is majorized by
\[
\begin{aligned} & \,\,\tilde{u}_{2}(\mathbf{x},\mathbf{x}^{(l)})\\
= & \,\,(\lambda_{u}-\lambda_{B})\mathbf{x}^{H}\mathbf{x}\\
 & +2{\rm Re}\left(\mathbf{x}^{H}\big(\mathbf{R}\!-\!\mathbf{B}\circ\big(\mathbf{x}^{(l)}(\mathbf{x}^{(l)})^{H}\big)\!-\!(\lambda_{u}\!-\!\lambda_{B})\mathbf{I}\big)\mathbf{x}^{(l)}\right)\\
 & +(\mathbf{x}^{(l)})^{H}((\lambda_{u}-\lambda_{B})\mathbf{I}-\mathbf{R}\!+\!\mathbf{B}\circ\big(\mathbf{x}^{(l)}(\mathbf{x}^{(l)})^{H}\big))\mathbf{x}^{(l)}.
\end{aligned}
\]
By ignoring the constant terms, the majorized problem of \eqref{eq:toep-diag}
is given by 
\begin{equation}
\begin{array}{ll}
\underset{\mathbf{x}}{\mathsf{minimize}} & {\rm Re}\left(\mathbf{x}^{H}\big(\mathbf{R}\!-\!\mathbf{B}\circ\big(\mathbf{x}^{(l)}(\mathbf{x}^{(l)})^{H}\big)\!-\!(\lambda_{u}\!-\!\lambda_{B})\mathbf{I}\big)\mathbf{x}^{(l)}\right)\\
\mathsf{subject\;to} & \left|x_{n}\right|=1,\,n=1,\ldots,N.
\end{array}\label{eq:linear_prob_diag}
\end{equation}
Similar to the problem \eqref{eq:linear_prob} in the previous section,
\eqref{eq:linear_prob_diag} also admits a closed form solution given
by 
\begin{equation}
x_{n}=e^{j{\rm arg}(\tilde{y}_{n})},\,n=1,\ldots,N,\label{eq:WISL_minimizer_diag}
\end{equation}
where
\begin{equation}
\begin{aligned}\tilde{\mathbf{y}} & =(\lambda_{u}\!-\!\lambda_{B})\mathbf{x}^{(l)}\!+\!\big(\mathbf{B}\circ\big(\mathbf{x}^{(l)}(\mathbf{x}^{(l)})^{H}\big)\big)\mathbf{x}^{(l)}\!-\!\mathbf{R}\mathbf{x}^{(l)}.\end{aligned}
\label{eq:y_tilde}
\end{equation}
It is worth noting that the $\tilde{\mathbf{y}}$ in \eqref{eq:y_tilde}
can be computed efficiently via FFT operations. To see that, we first
note that $\mathbf{R}\mathbf{x}^{(l)}$ can be computed by means of
FFT as described in the previous section. According to Lemma \ref{lem:Hadamard},
we have
\begin{equation}
\begin{aligned} & \,\,\big(\mathbf{B}\circ\big(\mathbf{x}^{(l)}(\mathbf{x}^{(l)})^{H}\big)\big)\mathbf{x}^{(l)}\\
= & \,\,{\rm diag}\big(\mathbf{B}{\rm Diag}(\mathbf{x}^{(l)})\big(\mathbf{x}^{(l)}(\mathbf{x}^{(l)})^{H}\big)^{T}\big)\\
= & \,\,{\rm diag}\big(\mathbf{B}\big(\mathbf{x}^{(l)}\circ\big(\mathbf{x}^{(l)}\big)^{*})\big(\mathbf{x}^{(l)}\big)^{T}\big)\\
= & \,\,{\rm diag}\big(\mathbf{B}\mathbf{1}\big(\mathbf{x}^{(l)}\big)^{T}\big)\\
= & \,\,(\mathbf{B}\mathbf{1})\circ\mathbf{x}^{(l)}.
\end{aligned}
\label{eq:Bx}
\end{equation}
Since $\mathbf{B}$ is symmetric Toeplitz, by Lemma \ref{lem:diagonal}
we can decompose $\mathbf{B}$ as 
\begin{equation}
\mathbf{B}=\frac{1}{2N}\mathbf{F}_{:,1:N}^{H}{\rm Diag}(\mathbf{F}\tilde{\mathbf{w}})\mathbf{F}_{:,1:N},
\end{equation}
and thus $\mathbf{B}\mathbf{1}$ can also be computed efficiently
via FFT operations. We further note that we only need to compute $\mathbf{B}\mathbf{1}$
once. The overall algorithm is then summarized in Algorithm \ref{alg:MWISL-Diag},
for which the main computation of each iteration is still just four
FFT (IFFT) operations and thus of order $\mathcal{O}(N\log N)$.

\begin{algorithm}[tbh]
\begin{algor}[1]
\item [{Require:}] \begin{raggedright}
sequence length $N$, weights $\{w_{k}\geq0\}_{k=1}^{N-1}$
\par\end{raggedright}
\item [{{*}}] Set $l=0$, initialize $\mathbf{x}^{(0)}$. 
\item [{{*}}] \begin{raggedright}
$\tilde{\mathbf{w}}=[0,w_{1}(N-1),\ldots,w_{N-1},0,w_{N-1},\ldots,w_{1}(N-1)]^{T}$
\par\end{raggedright}
\item [{{*}}] $\boldsymbol{\nu}=\mathbf{F}\tilde{\mathbf{w}}$
\item [{{*}}] $\mathbf{p}=\mathbf{F}_{:,1:N}^{H}\big(\boldsymbol{\nu}\circ\big(\mathbf{F}_{:,1:N}\mathbf{1}\big)\big)$
\item [{{*}}] $\lambda_{B}=\frac{1}{2}\left(\min_{1\leq i\leq N}\nu_{2i}+\min_{1\leq i\leq N}\nu_{2i-1}\right)$
\item [{repeat}]~

\begin{algor}[1]
\item [{{*}}] $\mathbf{f}=\mathbf{F}[\mathbf{x}^{(l)T},\mathbf{0}_{1\times N}]^{T}$
\item [{{*}}] $\mathbf{r}=\frac{1}{2N}\mathbf{F}^{H}\left|\mathbf{f}\right|^{2}$
\item [{{*}}] $\mathbf{c}=\mathbf{r}\circ[0,w_{1},\ldots,w_{N-1},0,w_{N-1},\ldots,w_{1}]^{T}$
\item [{{*}}] $\boldsymbol{\mu}=\mathbf{F}\mathbf{c}$
\item [{{*}}] $\lambda_{u}=\frac{1}{2}\big(\underset{1\leq i\leq N}{\max}\mu_{2i}+\underset{1\leq i\leq N}{\max}\mu_{2i-1}\big)$
\item [{{*}}] $\tilde{\mathbf{y}}=\mathbf{x}^{(l)}\!+\!\frac{\mathbf{p}\circ\mathbf{x}^{(l)}\!-\!\mathbf{F}_{:,1:N}^{H}(\boldsymbol{\mu}\circ\mathbf{f})}{2N(\lambda_{u}\!-\!\lambda_{B})}$
\item [{{*}}] $x_{n}^{(l+1)}=e^{j{\rm arg}(\tilde{y}_{n})},\,n=1,\ldots,N$ 
\item [{{*}}] $l\leftarrow l+1$ 
\end{algor}
\item [{until}] convergence
\end{algor}
\protect\caption{\label{alg:MWISL-Diag}MWISL-Diag - Monotonic minimizer for Weighted
ISL.}
\end{algorithm}

\section{Convergence Analysis and Acceleration Scheme\label{sec:Convergence-Acc}}

\subsection{Convergence Analysis}

The MWISL and MWISL-Diag algorithms given in Algorithm \ref{alg:MWISL}
and \ref{alg:MWISL-Diag} are based on the general majorization-minimization
framework, thus according to subsection \ref{sub:MM-Method}, we know
that the sequence of objective values (i.e., weighted ISL) evaluated
at $\{\mathbf{x}^{(l)}\}$ generated by the algorithms is nonincreasing.
And it is easy to see that the weighted ISL metric in \eqref{eq:WISL}
is bounded below by $0$, thus the sequence of objective values is
guaranteed to converge to a finite value. 

Now we further analyze the convergence property of the sequence $\{\mathbf{x}^{(l)}\}$
itself. In the following, we will focus on the sequence generated
by the MWISL algorithm (i.e, Algorithm \ref{alg:MWISL}), and prove
the convergence to a stationary point. The same result can be proved
for the MWISL-Diag algorithm (i.e, Algorithm \ref{alg:MWISL-Diag})
similarly.

To make it clear what is a stationary point in our case, we first
introduce a first-order optimality condition for minimizing a smooth
function over an arbitrary constraint set, which follows from \cite{Bertsekas2003}.
\begin{prop}
Let $f:\mathbf{R}^{n}\rightarrow\mathbf{R}$ be a smooth function,
and let $\mathbf{x}^{\star}$ be a local minimum of $f$ over a subset
$\mathcal{X}$ of $\mathbf{R}^{n}$. Then 
\begin{equation}
\nabla f(\mathbf{x}^{\star})^{T}\mathbf{z}\geq0,\,\forall\mathbf{z}\in T_{\mathcal{X}}(\mathbf{x}^{\star}),\label{eq:opt_condition}
\end{equation}
where $T_{\mathcal{X}}(\mathbf{x}^{\star})$ denotes the tangent cone
of $\mathcal{X}$ at $\mathbf{x}^{\star}.$
\end{prop}
A point $\mathbf{x}\in\mathcal{X}$ is said to be a stationary point
of the problem
\[
\begin{array}{ll}
\underset{\mathbf{x}\in\mathcal{X}}{\mathsf{minimize}} & f(\mathbf{x})\end{array}
\]
if it satisfies the first-order optimality condition \eqref{eq:opt_condition}.

To facilitate the analysis, we further note that upon defining 

\begin{equation}
\tilde{\mathbf{x}}=[{\rm Re}(\mathbf{x})^{T},{\rm Im}(\mathbf{x})^{T}]^{T},\label{eq:x_real}
\end{equation}
\begin{equation}
\tilde{\mathbf{U}}_{k}=\frac{1}{2}\left[\begin{array}{cc}
\mathbf{U}_{k}+\mathbf{U}_{k}^{T} & \mathbf{0}\\
\mathbf{0} & \mathbf{U}_{k}+\mathbf{U}_{k}^{T}
\end{array}\right],\label{eq:U_k_tilde}
\end{equation}
\begin{equation}
\hat{\mathbf{U}}_{k}=\frac{1}{2}\left[\begin{array}{cc}
\mathbf{0} & \mathbf{U}_{k}-\mathbf{U}_{k}^{T}\\
\mathbf{U}_{k}^{T}-\mathbf{U}_{k} & \mathbf{0}
\end{array}\right],\label{eq:U_k_hat}
\end{equation}
and based on the expression of $r_{k}$ in \eqref{eq:r_k}, it is
straightforward to show that the complex WISL minimization problem
\eqref{eq:WISL-minimize} is equivalent to the following real one:
\begin{equation}
\begin{array}{ll}
\underset{\tilde{\mathbf{x}}}{\mathsf{minimize}} & {\displaystyle \sum_{k=1}^{N-1}}w_{k}\left(\left(\tilde{\mathbf{x}}^{T}\tilde{\mathbf{U}}_{k}\tilde{\mathbf{x}}\right)^{2}+\left(\tilde{\mathbf{x}}^{T}\hat{\mathbf{U}}_{k}\tilde{\mathbf{x}}\right)^{2}\right)\\
\mathsf{subject\;to} & \tilde{x}_{n}^{2}+\tilde{x}_{n+N}^{2}=1,\,n=1,\ldots,N.
\end{array}\label{eq:prob_real_equivalent}
\end{equation}

We are now ready to state the convergence properties of MWISL.
\begin{thm}
\label{thm:MWISL-converge}Let $\{\mathbf{x}^{(l)}\}$ be the sequence
generated by the MWISL algorithm in Algorithm \ref{alg:MWISL}. Then
every limit point of the sequence $\{\mathbf{x}^{(l)}\}$ is a stationary
point of the problem \eqref{eq:WISL-minimize}.\end{thm}
\begin{IEEEproof}
Denote the objective functions of the problem \eqref{eq:WISL-minimize}
and its real equivalent \eqref{eq:prob_real_equivalent} by $f(\mathbf{x})$
and $\tilde{f}(\tilde{\mathbf{x}})$, respectively. Denote the constraint
sets of the problem \eqref{eq:WISL-minimize} and \eqref{eq:prob_real_equivalent}
by $\mathcal{C}$ and $\tilde{\mathcal{C}}$, respectively, i.e.,
$\mathcal{C}=\{\mathbf{x}\in\mathbf{C}^{N}|\left|x_{n}\right|=1,n=1,\ldots,N\}$
and $\tilde{\mathcal{C}}=\{\tilde{\mathbf{x}}\in\mathbf{R}^{2N}|\tilde{x}_{n}^{2}+\tilde{x}_{n+N}^{2}=1,\,n=1,\ldots,N\}$.
From the derivation of MWISL in subsection \ref{sub:WISL-Minimize_sub1},
we know that, at iteration $l$, $f(\mathbf{x})$ is majorized by
the function $u(\mathbf{x},\mathbf{x}^{(l)})$ in \eqref{eq:upper_u}
at $\mathbf{x}^{(l)}$ over $\mathcal{C}$. Then according to the
general MM scheme described in subsection \ref{sub:MM-Method}, we
have
\[
f(\mathbf{x}^{(l+1)})\leq u(\mathbf{x}^{(l+1)},\mathbf{x}^{(l)})\leq u(\mathbf{x}^{(l)},\mathbf{x}^{(l)})=f(\mathbf{x}^{(l)}),
\]
which means $\{f(\mathbf{x}^{(l)})\}$ is a nonincreasing sequence.

Since the sequence $\{\mathbf{x}^{(l)}\}$ is bounded, we know that
it has at least one limit point. Consider a limit point $\mathbf{x}^{(\infty)}$
and a subsequence $\{\mathbf{x}^{(l_{j})}\}$ that converges to $\mathbf{x}^{(\infty)},$
we have 
\[
\begin{aligned}u(\mathbf{x}^{(l_{j+1})},\mathbf{x}^{(l_{j+1})}) & =f(\mathbf{x}^{(l_{j+1})})\leq f(\mathbf{x}^{(l_{j}+1)})\\
 & \leq u(\mathbf{x}^{(l_{j}+1)},\mathbf{x}^{(l_{j})})\leq u(\mathbf{x},\mathbf{x}^{(l_{j})}),\forall\mathbf{x}\in\mathcal{C}.
\end{aligned}
\]
Letting $j\rightarrow+\infty,$ we obtain 
\begin{equation}
u(\mathbf{x}^{(\infty)},\mathbf{x}^{(\infty)})\leq u(\mathbf{x},\mathbf{x}^{(\infty)}),\,\forall\mathbf{x}\in\mathcal{C},
\end{equation}
i.e., $\mathbf{x}^{(\infty)}$ is a global minimizer of $u(\mathbf{x},\mathbf{x}^{(\infty)})$
over $\mathcal{C}$. With the definitions of $\tilde{\mathbf{x}}$,
$\tilde{\mathbf{U}}_{k}$ and $\hat{\mathbf{U}}_{k}$ given in \eqref{eq:x_real},
\eqref{eq:U_k_tilde} and \eqref{eq:U_k_hat}, and by ignoring the
constant terms in $u(\mathbf{x},\mathbf{x}^{(\infty)})$, it is easy
to show that minimizing $u(\mathbf{x},\mathbf{x}^{(\infty)})$ over
$\mathcal{C}$ is equivalent to the following real problem: 
\begin{equation}
\begin{array}{ll}
\underset{\tilde{\mathbf{x}}}{\mathsf{minimize}} & 4\tilde{\mathbf{x}}^{T}\mathbf{d}-2\big(\lambda_{{\rm max}}(\mathbf{L})N+\lambda_{u}\big)\tilde{\mathbf{x}}^{T}\tilde{\mathbf{x}}^{(\infty)}\\
\mathsf{subject\;to} & \tilde{\mathbf{x}}\in\tilde{\mathcal{C}},
\end{array}\label{eq:prob-mm-real}
\end{equation}
where $\tilde{\mathbf{x}}^{(\infty)}=[{\rm Re}(\mathbf{x}^{(\infty)})^{T},{\rm Im}(\mathbf{x}^{(\infty)})^{T}]^{T}$
and 
\begin{equation}
\mathbf{d}=\sum_{k=1}^{N-1}w_{k}\left(\tilde{\mathbf{x}}^{(\infty)T}\tilde{\mathbf{U}}_{k}\tilde{\mathbf{x}}^{(\infty)}\tilde{\mathbf{U}}_{k}+\tilde{\mathbf{x}}^{(\infty)T}\hat{\mathbf{U}}_{k}\tilde{\mathbf{x}}^{(\infty)}\hat{\mathbf{U}}_{k}\right)\tilde{\mathbf{x}}^{(\infty)}.
\end{equation}
Since $\mathbf{x}^{(\infty)}$ minimizes $u(\mathbf{x},\mathbf{x}^{(\infty)})$
over $\mathcal{C}$, $\tilde{\mathbf{x}}^{(\infty)}$ is a global
minimizer of \eqref{eq:prob-mm-real} and since $\tilde{\mathbf{x}}^{T}\tilde{\mathbf{x}}$
is just a constant over $\tilde{\mathcal{C}}$, $\tilde{\mathbf{x}}^{(\infty)}$
is also a global minimizer of 
\begin{equation}
\begin{array}{ll}
\underset{\tilde{\mathbf{x}}}{\mathsf{minimize}} & 4\tilde{\mathbf{x}}^{T}\mathbf{d}-\big(\lambda_{{\rm max}}(\mathbf{L})N+\lambda_{u}\big)\left(2\tilde{\mathbf{x}}^{T}\tilde{\mathbf{x}}^{(\infty)}-\tilde{\mathbf{x}}^{T}\tilde{\mathbf{x}}\right)\\
\mathsf{subject\;to} & \tilde{\mathbf{x}}\in\tilde{\mathcal{C}}.
\end{array}\label{eq:prob-real2}
\end{equation}
Then as a necessary condition, we have 
\begin{equation}
\nabla\tilde{u}(\tilde{\mathbf{x}}^{(\infty)})^{T}\mathbf{z}\geq0,\,\forall\mathbf{z}\in T_{\tilde{\mathcal{C}}}(\tilde{\mathbf{x}}^{(\infty)}),
\end{equation}
where $\tilde{u}(\tilde{\mathbf{x}})$ denotes the objective function
of \eqref{eq:prob-real2}. It is easy to check that 
\[
\nabla\tilde{f}(\tilde{\mathbf{x}}^{(\infty)})=\nabla\tilde{u}(\tilde{\mathbf{x}}^{(\infty)})=4\mathbf{d}.
\]
Thus we have 
\begin{equation}
\nabla\tilde{f}(\tilde{\mathbf{x}}^{(\infty)})^{T}\mathbf{z}\geq0,\,\forall\mathbf{z}\in T_{\tilde{\mathcal{C}}}(\tilde{\mathbf{x}}^{(\infty)}),
\end{equation}
implying that $\tilde{\mathbf{x}}^{(\infty)}$ is a stationary point
of the problem \eqref{eq:prob_real_equivalent}. Due to the equivalence
of problem \eqref{eq:prob_real_equivalent} and \eqref{eq:WISL-minimize},
the proof is complete.
\end{IEEEproof}

\subsection{Acceleration Scheme \label{sec:Acceleration-Schemes}}

In MM algorithms, the convergence speed is usually dictated by the
nature of the majorization functions. Due to the successive majorization
steps that we have carried out in the previous sections to construct
the majorization functions, the convergence of MWISL and MWISL-Diag
seems to be slow. In this subsection, we briefly introduce an acceleration
scheme that can be applied to accelerate the proposed MM algorithms.
It is the so called squared iterative method (SQUAREM), which was
originally proposed in \cite{SQUAREM} to accelerate any Expectation\textendash Maximization
(EM) algorithms. SQUAREM adapts the idea of the Cauchy-Barzilai-Borwein
(CBB) method \cite{raydan2002CBB}, which combines the classical steepest
descent method and the two-point step size gradient method \cite{barzilai1988_BB},
to solve the nonlinear fixed-point problem of EM. It only requires
the EM updating scheme and can be readily implemented as an off-the-shelf
accelerator. Since MM is a generalization of EM and the update rule
is also just a fixed-point iteration, SQUAREM can be easily applied
to MM algorithms after some minor modifications.

Suppose we have derived an MM algorithm to minimize $f(\mathbf{x})$
over $\mathcal{X}\subseteq\mathbf{C}^{n}$ and let $\mathbf{F}_{{\rm MM}}(\cdot)$
denote the nonlinear fixed-point iteration map of the MM algorithm:
\begin{equation}
\mathbf{x}^{(k+1)}=\mathbf{F_{{\rm MM}}}(\mathbf{x}^{(k)}).\label{eq:iter_map}
\end{equation}
For example, the iteration map of the MWISL algorithm is given by
\eqref{eq:WISL_minimizer}. Then the steps of the accelerated MM algorithm
based on SQUAREM are given in Algorithm \ref{alg:acc-MM}. A problem
of the general SQUAREM is that it may violate the nonlinear constraints,
so in Algorithm \ref{alg:acc-MM} we need to project wayward points
back to the feasible region by $\mbox{\ensuremath{\mathcal{P}}}_{\mathcal{X}}(\cdot)$.
For the unit-modulus constraints in the problem under consideration,
the projection can be done by simply applying the function $e^{j{\rm arg}(\cdot)}$
element-wise to the solution vectors. A second problem of SQUAREM
is that it can violate the descent property of the original MM algorithm.
To ensure the descent property, a strategy based on backtracking has
been adopted in Algorithm \ref{alg:acc-MM}, which repeatedly halves
the distance between $\alpha$ and $-1$:$\alpha\leftarrow(\alpha-1)/2$
until the descent property is maintained. To see why this works, we
first note that $\mbox{\ensuremath{\mathcal{P}}}_{\mathcal{X}}\left(\mathbf{x}^{(k)}-2\alpha\mathbf{r}+\alpha^{2}\mathbf{v}\right)=\mathbf{x}_{2}$
if $\alpha=-1.$ In addition, since $f(\mathbf{x}_{2})\leq f(\mathbf{x}^{(k)})$
due to the descent property of original MM steps, $f(\mathbf{x})\leq f(\mathbf{x}^{(k)})$
is guaranteed to hold as $\alpha\rightarrow-1.$ It is worth mentioning
that, in practice, usually only a few back-tracking steps are needed
to maintain the monotonicity of the algorithm.

\begin{algorithm}[tbh]
\begin{algor}[1]
\item [{Require:}] parameters
\item [{{*}}] Set $k=0$, initialize $\mathbf{x}^{(0)}$. 
\item [{repeat}]~

\begin{algor}[1]
\item [{{*}}] $\mathbf{x}_{1}=\mathbf{F_{{\rm MM}}}\left(\mathbf{x}^{(k)}\right)$
\item [{{*}}] $\mathbf{x}_{2}=\mathbf{F_{{\rm MM}}}\left(\mathbf{x}_{1}\right)$
\item [{{*}}] $\mathbf{r}=\mathbf{x}_{1}-\mathbf{x}^{(k)}$
\item [{{*}}] $\mathbf{v}=\mathbf{x}_{2}-\mathbf{x}_{1}-\mathbf{r}$
\item [{{*}}] Compute the step-length $\alpha=-\frac{\left\Vert \mathbf{r}\right\Vert }{\left\Vert \mathbf{v}\right\Vert }$
\item [{{*}}] $\mathbf{x}=\mbox{\ensuremath{\mathcal{P}}}_{\mathcal{X}}\left(\mathbf{x}^{(k)}-2\alpha\mathbf{r}+\alpha^{2}\mathbf{v}\right)$
\item [{while}] $f(\mathbf{x})>f(\mathbf{x}^{(k)})$
\item [{{*}}] $\alpha\leftarrow(\alpha-1)/2$
\item [{{*}}] $\mathbf{x}=\mbox{\ensuremath{\mathcal{P}}}_{\mathcal{X}}\left(\mathbf{x}^{(k)}-2\alpha\mathbf{r}+\alpha^{2}\mathbf{v}\right)$ 
\item [{endwhile}]~
\item [{{*}}] $\mathbf{x}^{(k+1)}=\mathbf{x}$
\item [{{*}}] $k\leftarrow k+1$ 
\end{algor}
\item [{until}] convergence
\end{algor}
\protect\caption{\label{alg:acc-MM}The acceleration scheme for MM algorithms.}
\end{algorithm}

\section{Minimizing the $\ell_{p}$-norm of Autocorrelation Sidelobes\label{sec:Minimizing-Lp-norm}}

In previous sections, we have developed algorithms to minimize the
weighted ISL metric of a unit-modulus sequence. It is clear that the
(unweighted) ISL metric is just the squared $\ell_{2}$-norm of the
autocorrelation sidelobes and in this section we would like to consider
the more general $\ell_{p}$-norm metric of the autocorrelation sidelobes
defined as 
\begin{equation}
\left({\displaystyle \sum_{k=1}^{N-1}}|r_{k}|^{p}\right)^{1/p}
\end{equation}
with $2\leq p<\infty$. The motivation is that by choosing different
$p$ values, we may get different metrics of particular interest.
For instance, by choosing $p\rightarrow+\infty$, the $\ell_{p}$-norm
metric tends to the $\ell_{\infty}$-norm of the autocorrelation sidelobes,
which is known as the peak sidelobe level (PSL). So it is well motivated
to consider the more general $\ell_{p}$-norm ($2\leq p<\infty$)
metric minimization problem
\begin{equation}
\begin{array}{ll}
\underset{x_{n}}{\mathsf{minimize}} & \left({\displaystyle \sum_{k=1}^{N-1}}|r_{k}|^{p}\right)^{1/p}\\
\mathsf{subject\;to} & \left|x_{n}\right|=1,\,n=1,\ldots,N,
\end{array}\label{eq:Lp-minimize-original}
\end{equation}
which is equivalent to
\begin{equation}
\begin{array}{ll}
\underset{x_{n}}{\mathsf{minimize}} & {\displaystyle \sum_{k=1}^{N-1}}|r_{k}|^{p}\\
\mathsf{subject\;to} & \left|x_{n}\right|=1,\,n=1,\ldots,N.
\end{array}\label{eq:Lp-minimize}
\end{equation}
It is easy to see that by choosing $p=2$, problem \eqref{eq:Lp-minimize}
becomes the ISL minimization problem. It is also worth noting that
we can easily incorporate weights in the objective (i.e., weighted
$\ell_{p}$-norm) as in \eqref{eq:WISL}.

To tackle the problem \eqref{eq:Lp-minimize} via majorization-minimization,
we need to construct a majorization function of the objective and
the idea is to majorize each $|r_{k}|^{p},k=1,\ldots,N-1$ by a quadratic
function of $|r_{k}|$. It is clear that when $p>2$, it is impossible
to construct a global quadratic majorization function of $\left|r_{k}\right|^{p}$.
However, we can still majorize it by a quadratic function locally
based on the following lemma.
\begin{lem}
\label{lem:p2_norm}Let $f(x)=x^{p}$ with $p\geq2$ and $x\in[0,t]$.
Then for any given $x_{0}\in[0,t)$, $f(x)$ is majorized at $x_{0}$
over the interval $[0,t]$ by the following quadratic function 
\begin{equation}
ax^{2}+(px_{0}^{p-1}-2ax_{0})x+ax_{0}^{2}-(p-1)x_{0}^{p},\label{eq:quad_major}
\end{equation}
where 
\begin{equation}
a=\frac{t^{p}-x_{0}^{p}-px_{0}^{p-1}(t-x_{0})}{(t-x_{0})^{2}}.\label{eq:a_val}
\end{equation}
\end{lem}
\begin{IEEEproof}
See Appendix \ref{sec:Proof-of-lemma_p2norm}.
\end{IEEEproof}
Given $\left|r_{k}^{(l)}\right|$ at iteration $l,$ according to
Lemma \ref{lem:p2_norm}, we know that $\left|r_{k}\right|^{p}$ ($p\geq2$)
is majorized at $\left|r_{k}^{(l)}\right|$ over $[0,t]$ by 
\begin{equation}
a_{k}|r_{k}|^{2}+b_{k}\left|r_{k}\right|+a_{k}\left|r_{k}^{(l)}\right|^{2}-(p-1)\left|r_{k}^{(l)}\right|^{p},
\end{equation}
where 
\begin{eqnarray}
a_{k} & = & \frac{t^{p}-\left|r_{k}^{(l)}\right|^{p}-p\left|r_{k}^{(l)}\right|^{p-1}(t-\left|r_{k}^{(l)}\right|)}{(t-\left|r_{k}^{(l)}\right|)^{2}},\label{eq:a_k}\\
b_{k} & = & p\left|r_{k}^{(l)}\right|^{p-1}-2a_{k}\left|r_{k}^{(l)}\right|.\label{eq:b_k}
\end{eqnarray}
Since the objective decreases at every iteration in the MM framework,
at the current iteration $l$, it is sufficient to majorize $\left|r_{k}\right|^{p}$
over the set on which the objective is smaller, i.e., $\sum_{k=1}^{N-1}\left|r_{k}\right|^{p}\leq\sum_{k=1}^{N-1}\left|r_{k}^{(l)}\right|^{p}$,
which implies $\left|r_{k}\right|\leq\left(\sum_{k=1}^{N-1}\left|r_{k}^{(l)}\right|^{p}\right)^{\frac{1}{p}}$.
Hence we can choose $t=\left(\sum_{k=1}^{N-1}\left|r_{k}^{(l)}\right|^{p}\right)^{\frac{1}{p}}$
in \eqref{eq:a_k}. Then the majorized problem of \eqref{eq:Lp-minimize}
in this case is given by (ignoring the constant terms)
\begin{equation}
\begin{array}{ll}
\underset{x_{n}}{\mathsf{minimize}} & {\displaystyle \sum_{k=1}^{N-1}}\left(a_{k}|r_{k}|^{2}+b_{k}\left|r_{k}\right|\right)\\
\mathsf{subject\;to} & \left|x_{n}\right|=1,\,n=1,\ldots,N.
\end{array}\label{eq:PSL_major1}
\end{equation}
We can see that the first term $\sum_{k=1}^{N-1}a_{k}|r_{k}|^{2}$
in the objective is just the weighted ISL metric with weights $w_{k}=a_{k}$,
and thus can be further majorized as in Section \ref{sub:WISL-Minimize_sub1}.
Following the steps in Section \ref{sub:WISL-Minimize_sub1} until
\eqref{eq:toep} (i.e., just the first majorization), we know that
it is majorized at $\mathbf{x}^{(l)}$ by (with constant terms ignored)
\begin{equation}
\mathbf{x}^{H}\left(\mathbf{R}-\lambda_{{\rm max}}(\mathbf{L})\mathbf{x}^{(l)}(\mathbf{x}^{(l)})^{H}\right)\mathbf{x},\label{eq:major_term1}
\end{equation}
where $\mathbf{R}$ and $\mathbf{L}$ are defined in \eqref{eq:R_mat}
and \eqref{eq:matrix_L} with $w_{k}=a_{k}.$ For the second term,
since it can be shown that $b_{k}\leq0$, we have 
\begin{eqnarray}
{\displaystyle \sum_{k=1}^{N-1}}b_{k}\left|r_{k}\right| & \leq & {\displaystyle \sum_{k=1}^{N-1}}b_{k}{\rm Re}\left\{ r_{k}^{*}\frac{r_{k}^{(l)}}{\left|r_{k}^{(l)}\right|}\right\} \\
 & = & {\displaystyle \sum_{k=1}^{N-1}}b_{k}{\rm Re}\left\{ {\rm Tr}(\mathbf{U}_{-k}\mathbf{x}\mathbf{x}^{H})\frac{r_{k}^{(l)}}{\left|r_{k}^{(l)}\right|}\right\} \\
 & = & {\rm Re}\left\{ \mathbf{x}^{H}\left({\displaystyle \sum_{k=1}^{N-1}}b_{k}\frac{r_{k}^{(l)}}{\left|r_{k}^{(l)}\right|}\mathbf{U}_{-k}\right)\mathbf{x}\right\} \\
 & = & \frac{1}{2}\mathbf{x}^{H}\left({\displaystyle \sum_{k=1-N}^{N-1}}b_{k}\frac{r_{k}^{(l)}}{\left|r_{k}^{(l)}\right|}\mathbf{U}_{-k}\right)\mathbf{x},\label{eq:major_term2}
\end{eqnarray}
where $b_{-k}=b_{k},k=1,\ldots,N-1$, $b_{0}=0$. By adding the two
majorization functions, i.e., \eqref{eq:major_term1} and \eqref{eq:major_term2},
and defining 
\begin{equation}
\hat{w}_{-k}=\hat{w}_{k}=a_{k}+\frac{b_{k}}{2\left|r_{k}^{(l)}\right|}=\frac{p}{2}\left|r_{k}^{(l)}\right|^{p-2},k=1,\ldots,N-1,\label{eq:w_k_hat}
\end{equation}
we have the majorized problem of \eqref{eq:PSL_major1} given by
\begin{equation}
\begin{array}{ll}
\underset{x_{n}}{\mathsf{minimize}} & \mathbf{x}^{H}\left(\tilde{\mathbf{R}}-\lambda_{{\rm max}}(\mathbf{L})\mathbf{x}^{(l)}(\mathbf{x}^{(l)})^{H}\right)\mathbf{x}\\
\mathsf{subject\;to} & \left|x_{n}\right|=1,\,n=1,\ldots,N,
\end{array}\label{eq:PSL_major2}
\end{equation}
where 
\begin{equation}
\tilde{\mathbf{R}}={\displaystyle \sum_{k=1-N}^{N-1}}\hat{w}_{k}r_{-k}^{(l)}\mathbf{U}_{k}.
\end{equation}
We can see that the problem \eqref{eq:PSL_major2} has the same form
as \eqref{eq:toep}, then by following similar steps as in section
\ref{sub:WISL-Minimize_sub1}, we can perform one more majorization
step and get the majorized problem
\begin{equation}
\begin{array}{ll}
\underset{\mathbf{x}}{\mathsf{minimize}} & \left\Vert \mathbf{x}-\mathbf{y}\right\Vert _{2}\\
\mathsf{subject\;to} & \left|x_{n}\right|=1,\,n=1,\ldots,N,
\end{array}\label{eq:prob_PSL_step}
\end{equation}
where 
\begin{equation}
\mathbf{y}=\left(\lambda_{{\rm max}}(\mathbf{L})N+\lambda_{u}\right)\mathbf{x}^{(l)}-\tilde{\mathbf{R}}\mathbf{x}^{(l)},
\end{equation}
and this time $\lambda_{{\rm max}}(\mathbf{L})$ should be computed
based on weights $a_{k}$ in \eqref{eq:a_k} and $\lambda_{u}$ is
based on the weights $\hat{w}_{k}$ in \eqref{eq:w_k_hat}. As in
previous cases, we only need to solve \eqref{eq:prob_PSL_step} in
closed form at every iteration and the overall algorithm is summarized
in Algorithm \ref{alg:MPSL}. Note that, to avoid numerical issue,
we have used the normalized $a_{k}$ and $\hat{w}_{k}$ (i.e., divided
by $t^{p}$) in Algorithm \ref{alg:MPSL}, which is equivalent to
divide the objective in \eqref{eq:PSL_major1} by $t^{p}$ during
the derivation. It is also worth noting that the algorithm can be
accelerated by the scheme described in subsection \ref{sec:Acceleration-Schemes}.

\begin{algorithm}[tbh]
\begin{algor}[1]
\item [{Require:}] \begin{raggedright}
sequence length $N$, parameter $p\geq2$
\par\end{raggedright}
\item [{{*}}] Set $l=0$, initialize $\mathbf{x}^{(0)}$. 
\item [{repeat}]~

\begin{algor}[1]
\item [{{*}}] $\mathbf{f}=\mathbf{F}[\mathbf{x}^{(l)T},\mathbf{0}_{1\times N}]^{T}$
\item [{{*}}] $\mathbf{r}=\frac{1}{2N}\mathbf{F}^{H}\left|\mathbf{f}\right|^{2}$
\item [{{*}}] \begin{raggedright}
$t=\left\Vert \mathbf{r}_{2:N}\right\Vert _{p}$
\par\end{raggedright}
\item [{{*}}] \begin{raggedright}
$a_{k}=\frac{1+(p-1)\big(\frac{\left|r_{k+1}\right|}{t}\big)^{p}-p\big(\frac{\left|r_{k+1}\right|}{t}\big)^{p-1}}{(t-\left|r_{k+1}\right|)^{2}},k=1,\ldots,N-1$
\par\end{raggedright}
\item [{{*}}] $\hat{w}_{k}=\frac{p}{2t^{2}}\big(\frac{\left|r_{k+1}\right|}{t}\big)^{p-2},k=1,\ldots,N-1$
\item [{{*}}] $\lambda_{L}=\max_{k}\{a_{k}(N-k)|k=1,\ldots,N-1\}$
\item [{{*}}] $\tilde{\mathbf{c}}=\mathbf{r}\circ[0,\hat{w}_{1},\ldots,\hat{w}_{N-1},0,\hat{w}_{N-1},\ldots,\hat{w}_{1}]^{T}$
\item [{{*}}] $\tilde{\boldsymbol{\mu}}=\mathbf{F}\tilde{\mathbf{c}}$
\item [{{*}}] $\lambda_{u}=\frac{1}{2}\big(\underset{1\leq i\leq N}{\max}\tilde{\mu}_{2i}+\underset{1\leq i\leq N}{\max}\tilde{\mu}_{2i-1}\big)$
\item [{{*}}] $\mathbf{y}=\mathbf{x}^{(l)}-\frac{\mathbf{F}_{:,1:N}^{H}(\tilde{\boldsymbol{\mu}}\circ\mathbf{f})}{2N(\lambda_{L}N+\lambda_{u})}$
\item [{{*}}] $x_{n}^{(l+1)}=e^{j{\rm arg}(y_{n})},\,n=1,\ldots,N$ 
\item [{{*}}] $l\leftarrow l+1$ 
\end{algor}
\item [{until}] convergence
\end{algor}
\protect\caption{\label{alg:MPSL}Monotonic minimizer for the $\ell_{p}$-metric of
autocorrelation sidelobes $(p\geq2)$.}
\end{algorithm}

\section{Numerical Experiments\label{sec:Numerical-Experiments}}

To compare the performance of the proposed MWISL algorithm and its
variants with existing algorithms and to show the potential of proposed
algorithms in designing sequences for various scenarios, we present
some experimental results in this section. All experiments were performed
on a PC with a 3.20GHz i5-3470 CPU and 8GB RAM.

\subsection{Weighted ISL Minimization}

In this subsection, we give an example of applying the proposed MWISL
and MWISL-Diag algorithms (with and without acceleration) to design
sequences with low correlation sidelobes only at required lags and
compare the performance with the WeCAN algorithm \cite{stoica2009new}.
The Matlab code of the benchmark algorithm, i.e., WeCAN, was downloaded
from the website\footnote{http://www.sal.ufl.edu/book/} of the book
\cite{he2012waveform}.

\begin{figure}[tbh]
\centering{}\includegraphics[width=0.95\columnwidth]{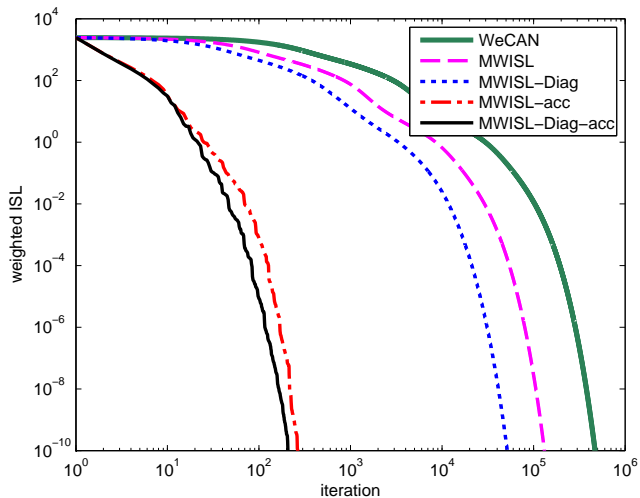}\protect\caption{\label{fig:WISL_vs_iteration}Evolution of the weighted ISL with respect
to the number of iterations.}
\end{figure}

\begin{figure}[tbh]
\centering{}\includegraphics[width=0.95\columnwidth]{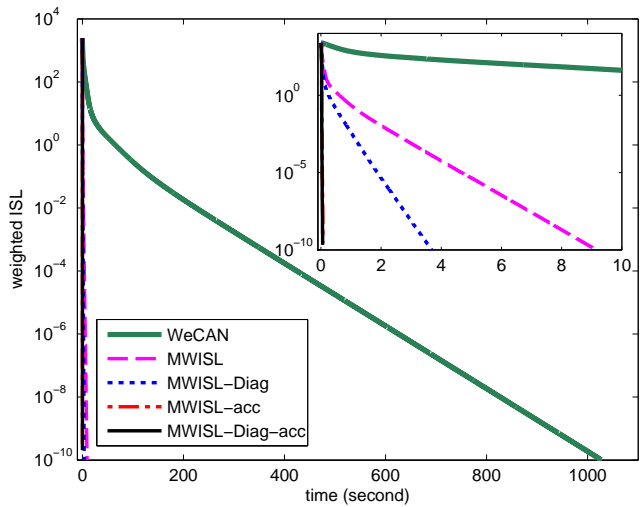}\protect\caption{\label{fig:WISL_vs_time}Evolution of the weighted ISL with respect
to time (in seconds). The plot within the time interval $[0,10]$
second is zoomed in and shown in the upper right corner.}
\end{figure}

\begin{figure}[htbp]
\centering{}\includegraphics[width=0.95\columnwidth]{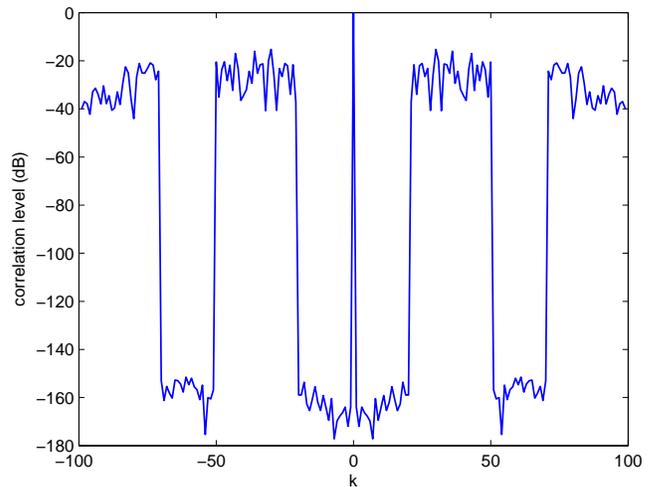}\protect\caption{\label{fig:WMISL_corr}Correlation level of the sequence of length
$N=100$ designed by accelerated MWISL-Diag algorithm with weights
in \eqref{eq:weights}.}
\end{figure}

Suppose we want to design a sequence of length $N=100$ and with small
correlations only at $r_{1},\ldots,r_{20}$ and $r_{51},\ldots,r_{70}.$
To achieve this, we set the weights $\{w_{k}\}_{k=1}^{99}$ as follows:
\begin{equation}
w_{k}=\begin{cases}
1, & k\in\{1,\dots,20\}\cup\{51,\ldots,70\}\\
0, & {\rm otherwise},
\end{cases}\label{eq:weights}
\end{equation}
such that only the autocorrelations at the required lags will be minimized.
Note that in this example there are $N-1=99$ degrees of freedom,
i.e., the free phases of $\{x_{n}\}_{n=1}^{N}$ (the initial phase
does not matter), and our goal is to match 80 real numbers (i.e.,
the real and imaginary parts of $r_{1},\ldots,r_{20}$ and $r_{51},\ldots,r_{70}$).
As there are enough degrees of freedom, the weighted ISL can be driven
to 0 in principle. So in this experiment, we will allow enough iterations
for the algorithms to be run and will not stop until the weighted
ISL goes below $10^{-10}$. The algorithms are initialized by a same
randomly generated sequence. The evolution curves of the weighted
ISL with respect to the number of iterations and time are shown in
Fig. \ref{fig:WISL_vs_iteration} and \ref{fig:WISL_vs_time}, respectively.
From Fig. \ref{fig:WISL_vs_iteration}, we can see that all the algorithms
can drive the weighted ISL to $10^{-10}$ when enough iterations are
allowed, but the proposed algorithms (especially the accelerated ones)
require far fewer iterations compared with the WeCAN algorithm. In
addition, we can see that the MWISL-Diag algorithms (with and without
acceleration) converge a bit faster than the corresponding MWISL algorithms,
which means the majorization function proposed in Section \ref{sec:Diag_majorize}
is somehow better than the one in Section \ref{sub:WISL-Minimize_sub1}.
From Fig. \ref{fig:WISL_vs_time}, we can see that in terms of the
computational time the superiority of the proposed algorithms is more
significant, more specially the accelerated MWISL and MWISL-Diag algorithms
take only 0.07 and 0.06 seconds respectively, while the WeCAN algorithm
takes more than 1000 seconds. It is because the proposed MWISL (and
MWISL-Diag) algorithms require only four FFT operations per iteration,
while each iteration of WeCAN requires $N$ computations of $2N$-point
FFTs. The correlation level of the output sequence of the accelerated
MWISL-Diag algorithm is shown in Fig. \ref{fig:WMISL_corr}, where
the correlation level is defined as 
\[
\textrm{correlation level}=20\log_{10}\left|\frac{r_{k}}{r_{0}}\right|,k=1-N,\ldots,N-1.
\]
We can see in Fig. \ref{fig:WMISL_corr} that the autocorrelation
sidelobes are suppressed to almost zero (about -160dB) at the required
lags.

\subsection{PSL Minimization}

\begin{figure}[t]
\centering{}\includegraphics[width=0.95\columnwidth]{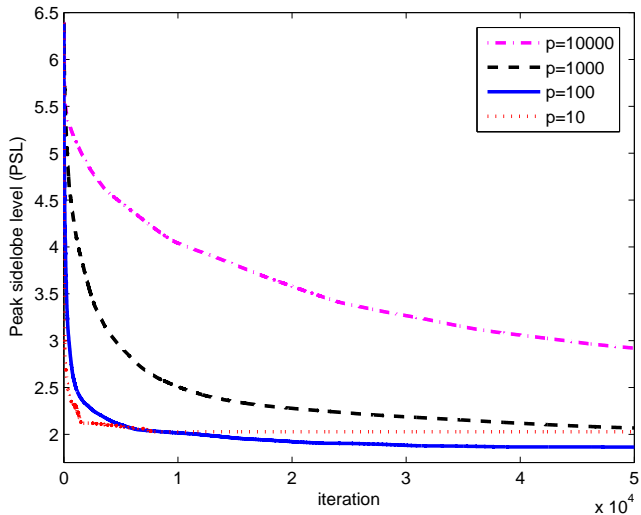}
\protect\caption{\label{fig:PSL_evolution}The evolution curves of the peak sidelobe
level (PSL).}
\end{figure}

In this subsection, we test the performance of the proposed Algorithm
\ref{alg:MPSL} in Section \ref{sec:Minimizing-Lp-norm} in minimizing
the peak sidelobe level (PSL) of the autocorrelation sidelobes, which
is of particular interest. To apply the algorihtm, we need to choose
the parameter $p.$ To examine the effect of the parameter $p,$ we
first apply the accelerated version of Algorithm \ref{alg:MPSL} (denoted
as MM-PSL) with four different $p$ values, i.e., $p=10,100,1000$
and $10000$, to design a sequence of length $N=400.$ Frank sequences
\cite{FrankSequence1963} are used to initialize the algorithm, which
are known to be sequences with good autocorrelation. More specifically,
Frank sequences are defined for lengths that are perfect squares and
the Frank sequence of length $N=M^{2}$ is given by 
\begin{equation}
x_{nM+k+1}=e^{j2\pi nk/M},n,k=0,1,\ldots,M-1.
\end{equation}
For all $p$ values, we stop the algorithm after $5\times10^{4}$
iterations and the evolution curves of the PSL are shown in Fig. \eqref{fig:PSL_evolution}.
From the figure, we can see that smaller $p$ values lead to faster
convergence. However, if $p$ is too small, it may not decrease the
PSL at a later stage, as we can see that $p=100$ finally gives smaller
PSL compared with $p=10.$ It may be explained by the fact that $\ell_{p}$-norm
with larger $p$ values approximates the $\ell_{\infty}$-norm better.
So in practice, gradually increasing the $p$ value is probably a
better approach.

In the second experiment, we consider both an increasing scheme of
$p$ (denoted as MM-PSL-adaptive) and the fixed $p$ scheme with $p=100$.
For the increasing scheme, we apply the MM-PSL algorithm with increasing
$p$ values $2,2^{2},\ldots,2^{13}$. For each $p$ value, the stopping
criterion was chosen to be $\left|{\rm obj}(\mathbf{x}^{(k+1)})-{\rm obj}(\mathbf{x}^{(k)})\right|/{\rm obj}(\mathbf{x}^{(k)})\leq10^{-5}/p$,
with ${\rm obj}(\mathbf{x})$ being the objective in \eqref{eq:Lp-minimize-original},
and the maximum allowed number of iterations was set to be $5\times10^{3}$.
For $p=2$, the algorithm is initialized by the Frank sequence and
for larger $p$ values, it is initialized by the solution obtained
at the previous $p$. For the fixed $p$ scheme, the stopping criterion
was chosen to be $\left|{\rm obj}(\mathbf{x}^{(k+1)})-{\rm obj}(\mathbf{x}^{(k)})\right|/{\rm obj}(\mathbf{x}^{(k)})\leq10^{-10}$,
and the maximum allowed number of iterations was $2\times10^{5}$.
In this case, in addition to the Frank sequence, the Golomb sequence
\cite{Golomb1993polyphase} was also used as the initial sequence,
which is also known for its good autocorrelation properties. In contrast
to Frank sequences, Golomb sequences are defined for any positive
integer and a Golomb sequence $\{x_{n}\}_{n=1}^{N}$ of length $N$
is given by\vspace{-10bp}

\begin{equation}
x_{n}=e^{j\pi(n-1)n/N},n=1,\ldots,N.
\end{equation}

\begin{figure}[tbph]
\centering{}\includegraphics[width=0.95\columnwidth]{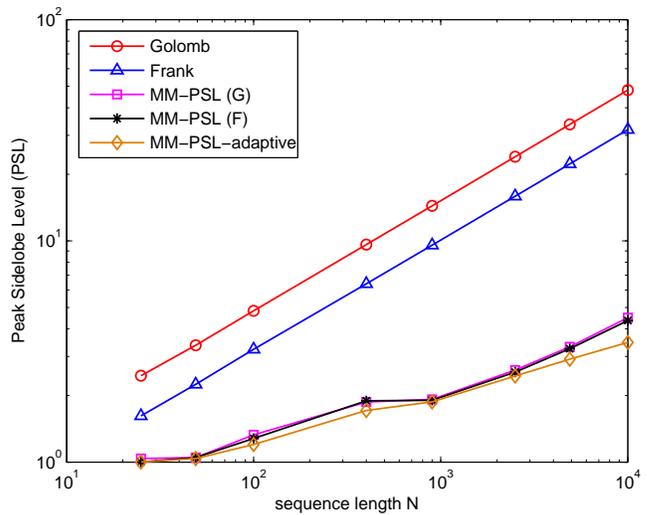}
\protect\caption{\label{fig:PSL}Peak sidelobe level (PSL) versus sequence length.
MM-PSL(G) and MM-PSL(F) denote the MM-PSL algorithm initialized by
Golomb and Frank sequences, respectively.}
\end{figure}

\begin{figure}[t]
\centering{}\includegraphics[width=0.95\columnwidth]{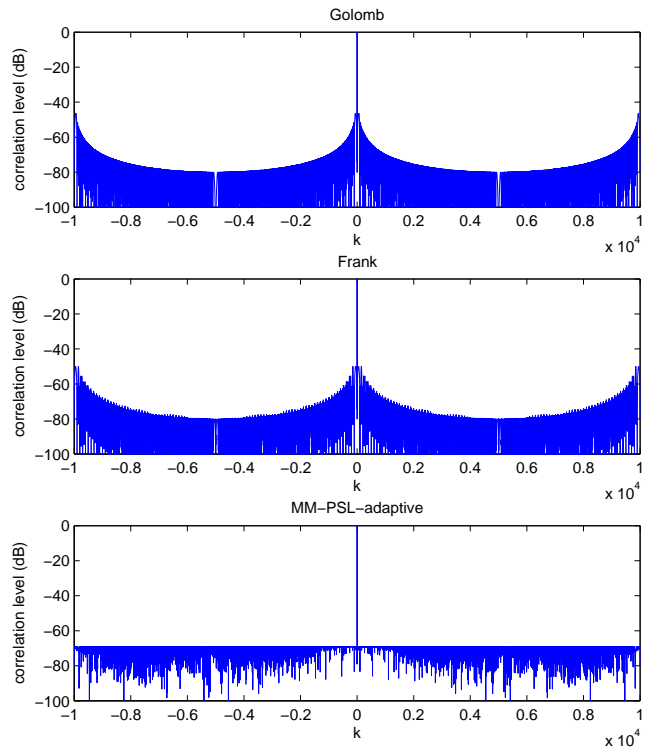}
\protect\caption{\label{fig:corr_level}Correlation level of the Golomb, Frank and
MM-PSL-adaptive sequences of length $N=10^{4}$.}
\vspace{-10bp}
\end{figure}

The two schemes are applied to design sequences of the following lengths:
$N=5^{2},7^{2},10^{2},20^{2},30^{2},50^{2},70^{2},100^{2}$, and the
PSL's of the resulting sequences are shown in Fig. \ref{fig:PSL}.
From the figure, we can see that for all lengths, the MM-PSL(G) and
MM-PSL(F) sequences give nearly the same PSL; both are much smaller
than the PSL of Golomb and Frank sequences, while a bit larger than
the PSL of MM-PSL-adaptive sequences. For example, when $N=10^{4},$
the PSL values of the MM-PSL(F) and MM-PSL-adaptive sequences are
4.36 and 3.48, while the PSL values of Golomb and Frank sequences
are 48.03 and 31.84, respectively. The correlation level of the Golomb,
Frank and the MM-PSL-adaptive sequences are shown in Fig. \ref{fig:corr_level}.
We can notice that the autocorrelation sidelobes of the Golomb and
Frank sequences are relatively large for $k$ close to $0$ and $N-1$,
while the MM-PSL-adaptive sequence has much more uniform autocorrelation
sidelobes across all lags. 

\section{Conclusion\label{sec:Conclusion}}

We have developed two efficient algorithms for the minimization of
the weighted integrated sidelobe level (WISL) metric of unit-modulus
sequences. The proposed algorithms are derived based on applying two
successive majorization steps and we have proved that they will converge
to a stationary point of the original WISL minimization problem. By
performing one more majorization step in the derivations, we have
extended the proposed algorithms to tackle the problem of minimizing
the $\ell_{p}$-norm of the autocorrelation sidelobes. All the algorithms
can be implemented by means of FFT operations and thus are computationally
very efficient in practice. An acceleration scheme that can be used
to further speed up the proposed algorithms has also been considered.
By some numerical examples, we have shown that the proposed WISL minimization
algorithms can generate sequences with virtually zero autocorrelation
sidelobes in some specified lag intervals with much lower computational
cost compared with the state-of-the-art. It has also been observed
that the proposed $\ell_{p}$-metric minimization algorithm can produce
long sequences with much more uniform autocorrelation sidelobes and
much smaller PSL compared with Frank and Golomb sequences, which are
known for their good autocorrelation properties.

\appendices{}

\section{Proof of Lemma \ref{lem:p2_norm} \label{sec:Proof-of-lemma_p2norm}}
\begin{IEEEproof}
For any given $x_{0}\in[0,t)$, let us consider a quadratic function
of the following form 
\begin{equation}
g(x|x_{0})=f(x_{0})+f^{\prime}(x_{0})(x-x_{0})+a(x-x_{0})^{2},\label{eq:g_x}
\end{equation}
where $a>0$. It is easy to check that $f(x_{0})=g(x_{0}|x_{0}).$
So to make $g(x|x_{0})$ be a majorization function of $f(x)$ at
$x_{0}$ over the interval $[0,t]$, we need to further have $f(x)\leq g(x|x_{0})$
for all $x\in[0,t],x\neq x_{0}$. Equivalently, we must have 
\begin{equation}
a\geq\frac{f(x)-f(x_{0})-f^{\prime}(x_{0})(x-x_{0})}{(x-x_{0})^{2}}
\end{equation}
for all $x\in[0,t]$, $x\neq x_{0}.$ Let us define the function 
\begin{equation}
A(x|x_{0})=\frac{f(x)-f(x_{0})-f^{\prime}(x_{0})(x-x_{0})}{(x-x_{0})^{2}}
\end{equation}
for all $x\neq x_{0}.$ The derivative of $A(x|x_{0})$ is given by
\[
A^{\prime}(x|x_{0})=\frac{f^{\prime}(x)+f^{\prime}(x_{0})-2(f(x)-f(x_{0}))/(x-x_{0})}{(x-x_{0})^{2}}.
\]
Since $f^{\prime}(x)=px^{p-1}$ is convex on $[0,t]$ when $p\geq2$,
we have
\begin{eqnarray*}
\frac{f(x)-f(x_{0})}{x-x_{0}} & = & \int_{0}^{1}f^{\prime}(x_{0}+\tau(x-x_{0}))d\tau\\
 & \leq & \int_{0}^{1}\left(f^{\prime}(x_{0})+\tau\left(f^{\prime}(x)-f^{\prime}(x_{0})\right)\right)d\tau\\
 & = & \frac{1}{2}\left(f^{\prime}(x)+f^{\prime}(x_{0})\right),
\end{eqnarray*}
which implies $A^{\prime}(x|x_{0})\geq0$ for all $x\in[0,t],x\neq x_{0}.$
Thus, $A(x|x_{0})$ is increasing on the interval $[0,t]$ and the
maximum is achieved at $x=t.$ Then the smallest $a$ we may choose
is
\begin{equation}
\begin{aligned}a & =\max_{x\in[0,t],x\neq x_{0}}\,A(x|x_{0})\\
 & =\frac{t^{p}-x_{0}^{p}-px_{0}^{p-1}(t-x_{0})}{(t-x_{0})^{2}}.
\end{aligned}
\end{equation}
By substituting $a$ into $g(x|x_{0})$ in \eqref{eq:g_x} and appropriately
rearranging terms, we can obtain the function in \eqref{eq:quad_major}.
\end{IEEEproof}
\vspace{-5bp}
\bibliographystyle{IEEEtran}
\bibliography{WISL}

\begin{thebibliography}{10}
\providecommand{\url}[1]{#1}
\csname url@samestyle\endcsname
\providecommand{\newblock}{\relax}
\providecommand{\bibinfo}[2]{#2}
\providecommand{\BIBentrySTDinterwordspacing}{\spaceskip=0pt\relax}
\providecommand{\BIBentryALTinterwordstretchfactor}{4}
\providecommand{\BIBentryALTinterwordspacing}{\spaceskip=\fontdimen2\font plus
\BIBentryALTinterwordstretchfactor\fontdimen3\font minus
  \fontdimen4\font\relax}
\providecommand{\BIBforeignlanguage}[2]{{%
\expandafter\ifx\csname l@#1\endcsname\relax
\typeout{** WARNING: IEEEtran.bst: No hyphenation pattern has been}%
\typeout{** loaded for the language `#1'. Using the pattern for}%
\typeout{** the default language instead.}%
\else
\language=\csname l@#1\endcsname
\fi
#2}}
\providecommand{\BIBdecl}{\relax}
\BIBdecl

\bibitem{turyn1968sequences}
R.~Turyn, ``Sequences with small correlation,'' \emph{Error correcting codes},
  pp. 195--228, 1968.

\bibitem{channel_estimation}
P.~Spasojevic and C.~Georghiades, ``Complementary sequences for {ISI} channel
  estimation,'' \emph{IEEE Transactions on Information Theory}, vol.~47, no.~3,
  pp. 1145--1152, Mar. 2001.

\bibitem{levanon2004radar}
N.~Levanon and E.~Mozeson, \emph{{Radar Signals}}.\hskip 1em plus 0.5em minus
  0.4em\relax John Wiley \& Sons, 2004.

\bibitem{golomb2005signal}
S.~W. Golomb and G.~Gong, \emph{{Signal Design for Good Correlation: For
  Wireless Communication, Cryptography, and Radar}}.\hskip 1em plus 0.5em minus
  0.4em\relax Cambridge University Press, 2005.

\bibitem{he2012waveform}
H.~He, J.~Li, and P.~Stoica, \emph{{Waveform Design for Active Sensing Systems:
  A Computational Approach}}.\hskip 1em plus 0.5em minus 0.4em\relax Cambridge
  University Press, 2012.

\bibitem{barker1953group}
R.~Barker, ``Group synchronizing of binary digital systems,''
  \emph{Communication theory}, pp. 273--287, 1953.

\bibitem{golomb1965generalized}
S.~W. Golomb and R.~A. Scholtz, ``Generalized {Barker} sequences,'' \emph{IEEE
  Transactions on Information Theory}, vol.~11, no.~4, pp. 533--537, 1965.

\bibitem{zhang1989polyBarker}
N.~Zhang and S.~W. Golomb, ``Sixty-phase generalized {Barker} sequences,''
  \emph{IEEE Transactions on Information Theory}, vol.~35, no.~4, pp. 911--912,
  1989.

\bibitem{friese1994polyBarker}
M.~Friese and H.~Zottmann, ``Polyphase {Barker} sequences up to length 31,''
  \emph{Electronics letters}, vol.~30, no.~23, pp. 1930--1931, 1994.

\bibitem{brenner1998polyphase}
A.~Brenner, ``Polyphase {Barker} sequences up to length 45 with small
  alphabets,'' \emph{Electronics letters}, vol.~34, no.~16, pp. 1576--1577,
  1998.

\bibitem{borwein2005polyphase}
P.~Borwein and R.~Ferguson, ``Polyphase sequences with low autocorrelation,''
  \emph{IEEE Transactions on Information Theory}, vol.~51, no.~4, pp.
  1564--1567, 2005.

\bibitem{Polycode2009}
C.~Nunn and G.~Coxson, ``Polyphase pulse compression codes with optimal peak
  and integrated sidelobes,'' \emph{IEEE Transactions on Aerospace and
  Electronic Systems}, vol.~45, no.~2, pp. 775--781, Apr. 2009.

\bibitem{FrankSequence1963}
R.~L. Frank, ``Polyphase codes with good nonperiodic correlation properties,''
  \emph{IEEE Transactions on Information Theory}, vol.~9, no.~1, pp. 43--45,
  Jan. 1963.

\bibitem{chu1972polyphase}
D.~Chu, ``Polyphase codes with good periodic correlation properties
  (corresp.),'' \emph{IEEE Transactions on Information Theory}, vol.~18, no.~4,
  pp. 531--532, 1972.

\bibitem{Golomb1993polyphase}
N.~Zhang and S.~W. Golomb, ``Polyphase sequence with low autocorrelations,''
  \emph{IEEE Transactions on Information Theory}, vol.~39, no.~3, pp.
  1085--1089, 1993.

\bibitem{turyn1967correlation}
R.~Turyn, ``The correlation function of a sequence of roots of 1 (corresp.),''
  \emph{IEEE Transactions on Information Theory}, vol.~13, no.~3, pp. 524--525,
  1967.

\bibitem{mow1997aperiodic}
W.~H. Mow and S.-Y. Li, ``Aperiodic autocorrelation and crosscorrelation of
  polyphase sequences,'' \emph{IEEE Transactions on Information Theory},
  vol.~43, no.~3, pp. 1000--1007, May 1997.

\bibitem{stoica2009new}
P.~Stoica, H.~He, and J.~Li, ``New algorithms for designing unimodular
  sequences with good correlation properties,'' \emph{IEEE Transactions on
  Signal Processing}, vol.~57, no.~4, pp. 1415--1425, 2009.

\bibitem{ITROX2012}
M.~Soltanalian and P.~Stoica, ``Computational design of sequences with good
  correlation properties,'' \emph{IEEE Transactions on Signal Processing},
  vol.~60, no.~5, pp. 2180--2193, May 2012.

\bibitem{naghsh2013unified}
M.~M. Naghsh, M.~Modarres-Hashemi, S.~ShahbazPanahi, M.~Soltanalian, and
  P.~Stoica, ``Unified optimization framework for multi-static radar code
  design using information-theoretic criteria,'' \emph{IEEE Transactions on
  Signal Processing}, vol.~61, no.~21, pp. 5401--5416, 2013.

\bibitem{UQP}
M.~Soltanalian and P.~Stoica, ``Designing unimodular codes via quadratic
  optimization,'' \emph{IEEE Transactions on Signal Processing}, vol.~62,
  no.~5, pp. 1221--1234, Mar. 2014.

\bibitem{MISL}
\BIBentryALTinterwordspacing
J.~Song, P.~Babu, and D.~Palomar, ``Optimization methods for designing
  sequences with low autocorrelation sidelobes,'' \emph{accepted in IEEE
  Transactions on Signal Processing}, to appear 2015. [Online]. Available:
  \url{http://arxiv.org/abs/1501.02252}
\BIBentrySTDinterwordspacing

\bibitem{torii2004new}
H.~Torii, M.~Nakamura, and N.~Suehiro, ``A new class of zero-correlation zone
  sequences,'' \emph{IEEE Transactions on Information Theory}, vol.~50, no.~3,
  pp. 559--565, 2004.

\bibitem{tang2008new}
X.~Tang and W.~H. Mow, ``A new systematic construction of zero correlation zone
  sequences based on interleaved perfect sequences,'' \emph{IEEE Transactions
  on Information Theory}, vol.~54, no.~12, pp. 5729--5734, 2008.

\bibitem{hunter2004MMtutorial}
D.~R. Hunter and K.~Lange, ``A tutorial on {MM} algorithms,'' \emph{The
  American Statistician}, vol.~58, no.~1, pp. 30--37, 2004.

\bibitem{MM_Stoica}
P.~Stoica and Y.~Selen, ``Cyclic minimizers, majorization techniques, and the
  expectation-maximization algorithm: a refresher,'' \emph{IEEE Signal
  Processing Magazine}, vol.~21, no.~1, pp. 112--114, Jan. 2004.

\bibitem{razaviyayn2013unified}
M.~Razaviyayn, M.~Hong, and Z.-Q. Luo, ``A unified convergence analysis of
  block successive minimization methods for nonsmooth optimization,''
  \emph{SIAM Journal on Optimization}, vol.~23, no.~2, pp. 1126--1153, 2013.

\bibitem{Aldo2013}
G.~Scutari, F.~Facchinei, P.~Song, D.~P. Palomar, and J.-S. Pang,
  ``Decomposition by partial linearization: Parallel optimization of
  multi-agent systems,'' \emph{IEEE Transactions on Signal Processing},
  vol.~62, no.~3, pp. 641--656, Feb. 2014.

\bibitem{eig_localization}
P.~J. S.~G. Ferreira, ``Localization of the eigenvalues of {Toeplitz} matrices
  using additive decomposition, embedding in circulants, and the {Fourier}
  transform,'' in \emph{Proceedings of the 10th IFAC Symposium on System
  Identification}, Copenhagen, Denmark, Jul. 1994, pp. 271--276.

\bibitem{Gray2006}
R.~M. Gray, \emph{{Toeplitz and Circulant Matrices: A review}}.\hskip 1em plus
  0.5em minus 0.4em\relax Now Publishers Inc, Jan. 2006, vol.~2, no.~3.

\bibitem{roger1994topics}
H.~Roger and R.~J. Charles, \emph{{Topics in Matrix Analysis}}.\hskip 1em plus
  0.5em minus 0.4em\relax Cambridge University Press, 1994.

\bibitem{Bertsekas2003}
D.~P. Bertsekas, A.~Nedi\'{c}, and A.~E. Ozdaglar, \emph{{Convex Analysis and
  Optimization}}.\hskip 1em plus 0.5em minus 0.4em\relax Athena Scientific,
  2003.

\bibitem{SQUAREM}
R.~Varadhan and C.~Roland, ``Simple and globally convergent methods for
  accelerating the convergence of any {EM} algorithm,'' \emph{Scandinavian
  Journal of Statistics}, vol.~35, no.~2, pp. 335--353, 2008.

\bibitem{raydan2002CBB}
M.~Raydan and B.~F. Svaiter, ``Relaxed steepest descent and
  {Cauchy-Barzilai-Borwein} method,'' \emph{Computational Optimization and
  Applications}, vol.~21, no.~2, pp. 155--167, 2002.

\bibitem{barzilai1988_BB}
J.~Barzilai and J.~M. Borwein, ``Two-point step size gradient methods,''
  \emph{IMA Journal of Numerical Analysis}, vol.~8, no.~1, pp. 141--148, 1988.

\end{thebibliography}

\end{document}